\documentclass[11pt,oneside,a4paper,reqno]{amsart}
\RequirePackage{fix-cm} 
\usepackage[T1]{fontenc}
\usepackage[english]{babel}
\usepackage{amssymb,amsmath,amsthm,graphicx,amsfonts}
\usepackage[foot]{amsaddr}
\usepackage{lmodern}
\usepackage[dvipsnames]{xcolor}
\definecolor{ourcolor}{RGB}{0,90,175}

\usepackage[
   rm={oldstyle=false,proportional=true},
   sf={oldstyle=true,proportional=true},
   tt={oldstyle=true,proportional=true,variable=true},
   qt=false,
 ]{cfr-lm}
\usepackage{cite,enumerate,setspace}
\usepackage[normalem]{ulem}
\usepackage[colorlinks=true,linkcolor=blue, 
            citecolor=blue, urlcolor=blue]{hyperref}
\usepackage{enumitem}

\numberwithin{equation}{section}

\usepackage[left=2.8cm,right=2.8cm,top=2.8cm,bottom=2.8cm]{geometry}
\let\originalleft\left
\let\originalright\right
\renewcommand{\left}{\mathopen{}\mathclose\bgroup\originalleft}
\renewcommand{\right}{\aftergroup\egroup\originalright}

 \makeatletter
 \def\@seccntformat#1{\hspace*{0mm}%
  \protect\textup{\protect\@secnumfont
    \ifnum\pdfstrcmp{subsection}{#1}=0 \bfseries\fi
    \csname the#1\endcsname
    \protect\@secnumpunct
      }%
 }

\newcommand{\assign}{:=}
\newcommand{\backassign}{=:}
\newcommand{\mathD}{\mathrm{D}}

\newcommand{\of}{:}
\newcommand{\tmSep}{; }
\newcommand{\tmabbr}[1]{#1}

\newcommand{\tmem}[1]{{\em #1\/}}
\newcommand{\tmname}[1]{\textup{\textsc{#1}}}
\newcommand{\tmop}[1]{\ensuremath{\operatorname{#1}}}
\newcommand{\tmsep}{, }

\newcommand{\tmtextit}[1]{\text{{\itshape{#1}}}}
\newcommand{\tmtextup}[1]{\text{{\upshape{#1}}}}

\newtheorem{corollary}{Corollary}[section]

{\theoremstyle{remark}\newtheorem{remark}{Remark}[section]}
\newtheorem{theorem}{Theorem}[section]

\newcommand{\de}{\mathrm{d}}
\newcommand{\di}{d_i}
\newcommand{\R}{\mathbb{R}}
\newcommand{\Stwo}{\mathbb{S}}
\newcommand{\ep}{\varepsilon}
\newcommand{\hyperrefstar}[2]{{\color{ourcolor}\tmtextup{#1}}}
\newcommand{\opt}{\tmtextup{(}}
\newcommand{\cpt}{\tmtextup{)}}
\newcommand\UUU{\color{black}}
\newcommand\EEE{\color{black}}

\title[A BBM formula accounting for antisymmetric exchange interactions]{A Bourgain--Brezis--Mironescu formula accounting for nonlocal antisymmetric exchange interactions}

\author[E. Davoli] {Elisa Davoli} 
\author[G. D. Fratta] {Giovanni Di Fratta} 
\author[R. Giorgio]{Rossella Giorgio}

\date{}
\AtEndDocument{
	\hrule
	\medskip
	\small{
		\noindent
		{\sc E. Davoli},
		Institute of Analysis and Scientific Computing,
		TU Wien,
		Wiedner Hauptstra\ss e 8-10, 1040 Vienna, Austria.
		
		\noindent
		E-mail:
		\href{mailto:elisa.davoli@tuwien.ac.at}{\tt elisa.davoli@tuwien.ac.at}.
	}

	\medskip \medskip \medskip \hrule  \medskip
	\small{
		\noindent
		{\sc G. Di Fratta},
		Dipartimento di Matematica e Applicazioni ``R. Caccioppoli'',
		Universit\`a degli Studi di Napoli ``Federico II'',
		Via Cintia, Complesso Monte S. Angelo, 80126 Naples, Italy.
		
		\noindent
		E-mail:
		\href{mailto:giovanni.difratta@unina.it}{\tt giovanni.difratta@unina.it}.
	}
	
	\medskip \medskip \medskip \hrule  \medskip
    \small{
		\noindent
		{\sc R. Giorgio},
		Institute of Analysis and Scientific Computing,
		TU Wien,
		Wiedner Hauptstra\ss e 8-10, 1040 Vienna, Austria
           {\sc and}
            \noindent MedUni Wien, Währinger Gürtel 18-20, 1090 Vienna, Austria.
            
		\noindent
		E-mail:
		\href{mailto:rossella.giorgio@student.tuwien.ac.at}{\tt rossella.giorgio@student.tuwien.ac.at}.
	}
	\medskip
\hrule
}

\begin{document}

\vskip .2truecm

\begin{abstract}
   The present study concerns the nonlocal-to-local convergence of a family of exchange energy functionals in the limit of very short-range interactions. The analysis accounts for both symmetric and antisymmetric exchange. Our result is twofold. First, we extend the Bourgain-Brezis-Mironescu formula to encompass the scenario where antisymmetric contributions are encoded into the energy. Second, we prove that, under physically relevant assumptions on the families of exchange kernels, the family of nonlocal functionals Gamma-converges to their local counterparts. As a byproduct of our analysis, we obtain a rigorous justification of Dzyaloshinskii--Moriya interactions in chiral magnets under the form commonly adopted in the variational theory of micromagnetism when modeling antisymmetric exchange interactions.
\end{abstract}

\subjclass{46E35\tmSep  49J45\tmSep  49S05}

\keywords{Nonlocal energies\tmsep  Bourgain-Brezis-Mironescu formula\tmsep 
$\Gamma$-Convergence\tmsep antisymmetric exchange interactions\tmsep Micromagnetics\tmsep Dzyaloshinskii--Moriya interaction (DMI)\tmsep  Magnetic skyrmions.}

{\maketitle}


\section{Introduction and motivation}
The present paper investigates the short-range interaction limit of a family
of nonlocal exchange energies of the form
\begin{equation}
  \mathcal{E}_{\ep} (m) \assign \mathcal{F}_{\ep} (m) +\mathcal{H}_{\ep} (m),
  \label{func:nonloc}
\end{equation}
\UUU where \EEE the energy functionals $\mathcal{F}_{\ep}$ and $\mathcal{H}_{\ep}$ \UUU are \EEE given
by
\begin{align}
     \mathcal{F}_{\ep} (m) & \assign \iint_{\Omega \times \Omega} \rho_{\ep} 
     (x - y)  \frac{|m (x) - m (y) |^2}{|x - y|^2} \de x \hspace{0.17em} \de
     y,  \label{eq:Fepssymexch}\\
     \mathcal{H}_{\ep} (m) & \assign \iint_{\Omega \times \Omega} \nu_{\ep} 
     (x - y) \cdot \frac{(m (x) \times m (y))}{|x - y|} \hspace{0.17em} \de x
     \hspace{0.17em} \de y,  \label{eq:Hepsasymexch}
\end{align}
\UUU and are both \EEE defined on a suitable (metric) subspace of $L^2 (\Omega ; \Stwo^2)$. Here,
$L^2 (\Omega ; \Stwo^2)$ denotes the space of vector-valued maps $m : \Omega
\to \Stwo^2$, where $\Omega$ is a bounded Lipschitz domain of $\R^3$ and
$\Stwo^2$ the unit sphere of $\R^3$.

The energy functional $\mathcal{E}_{\ep}$ is the result of two different types
of interactions. The term $\mathcal{F}_{\ep}$ in \eqref{eq:Fepssymexch}
accounts for the so-called \tmtextit{symmetric} exchange interactions, whereas
the energy term $\mathcal{H}_{\ep}$ in \eqref{eq:Hepsasymexch} accounts for
\tmtextit{antisymmetric} exchange interactions. The scalar kernel $\rho_{\ep}
: \R^3 \to \R_+$ and the vector-valued kernel $\nu_{\ep} : \R^3 \to \R^3$
model the strength and positional configuration of the exchange interactions
at spatial scale $\ep > 0$. They will be referred to as symmetric and
antisymmetric exchange kernels, respectively.

The main aim of this paper is to show that, under physically relevant
assumptions on the families of exchange kernels $\left( \rho_{\ep}
\right)_{\varepsilon}$ and $(\nu_{\varepsilon})_{\varepsilon}$, the family
$\mathcal{E}_{\ep}$ converges, in the sense of $\Gamma$-convergence in $L^2
(\Omega ; \Stwo^2)$, and up to constant factors, to the local energy
functional
\begin{equation}
  \mathcal{E} (m) \assign \mathcal{F} (m) +\mathcal{H} (m) \assign \frac{1}{2}
  \int_{\Omega} | \nabla m (x) |^2 + \sum_{i = 1}^3 \int_{\Omega} m (x) \cdot
  (d_i \times \partial_i m (x))  \de x, \label{func:loc0}
\end{equation}
provided that $m \in H^1 (\Omega ; \Stwo^2)$. Here, the quantities $d_i \in
\R^3$ denote constant vectors referred to as {\tmem{Dzyaloshinskii vectors}}
in the \UUU Physics \EEE literature; their expression strongly depends on the \UUU limiting \EEE
behavior of the family $(\nu_{\varepsilon})_{\varepsilon}$.

Rigorous statements will be given in Section~\ref{sec:rigsettingmainresults},
where we formulate explicit assumptions on the exchange kernels.

\subsection{Outline}The paper is organized as follows. In the rest of this
Section, we present a brief overview of the physical framework that motivated
our investigation and elucidate the significance of our findings in relation
to prior research on the subject. In Section~\ref{sec:rigsettingmainresults},
we present the precise formulation of the problem and provide a detailed
account of the contributions of our work. Proofs of our main results are given
in Sections~\ref{Sec:Point} and ~\ref{sec:Gammaconv}. Specifically,
Section~\ref{Sec:Point} focuses on the pointwise convergence of the
antisymmetric exchange interactions \eqref{eq:Hepsasymexch}, while
Section~\ref{sec:Gammaconv} establishes the $\Gamma$-convergence of the
nonlocal family $\left( \mathcal{E}_{\ep} \right)$ in \eqref{func:nonloc} to
the local energy functional $\mathcal{E}$ in \eqref{func:loc0}.

\subsection{Physics context: symmetric and antisymmetric
exchange interactions in micromagnetics}Reliable theoretical models for
studying magnetic phenomena must depend on the relevant length scales. At the
mesoscopic scale, there is a well-established and effective variational theory
of micromagnetism, whose roots may be found in the works of
Landau--Lifshitz~{\cite{LandauA1935}} and Brown~{\cite{BrownB1963,BrownB1962}}
on fine ferromagnetic particles. In this theory, for a rigid ferromagnetic
particle occupying a region $\Omega \subseteq \R^3$, the order parameter is
the magnetization field $M : \Omega \rightarrow \mathbb{} \R^3$. The modulus
of $M$, $M_s \assign | M |$, is called {\tmem{spontaneous magnetization}} and
is a function of the temperature that vanishes above the so-called
{\tmem{Curie point}} $T_c$: a critical value strongly depends on the specific
crystal structure of the ferromagnet. When the specimen is at a fixed
temperature well below $T_c$, the function $M_s$ can be assumed constant in
$\Omega$, and the magnetization can be conveniently written as $M \assign M_s
m$, where $m : \Omega \rightarrow \Stwo^2$ is a vector field with values in
the unit sphere of $\R^3$ (cf.~{\cite{BrownB1963,hubert2008magnetic}}).

Despite $| M |$ being constant in $\Omega$, this is generally not the case for
the direction of $M$, and according to the variational theory of
micromagnetism, the observable magnetization patterns are the local minimizers
of the {\tmem{micromagnetic energy functional}}, which, after normalization,
reads as\footnote{In writing {\eqref{eq:GLunorm}{{\tmem{ is non-convex,
non-local, and contains multiple length scales}}}}, we neglected the
{\tmem{magnetocrystalline anisotropy}} and {\tmem{Zeeman}} energy, but only to
shorten the notation. Indeed, although these contributions are of fundamental
importance in ferromagnetism~{\cite{BrownB1963,BrownB1962}}, from the
variational point of view, they behave like continuous perturbations, their
analysis is usually straightforward, and in our specific context, they play no
role.}
\begin{equation}
  \mathcal{G} (m) \assign \mathcal{F} (m) +\mathcal{W} (m) \assign \frac{1}{2}
  \int_{\Omega} | \nabla m|^2 + \frac{1}{2}  \int_{\R^3} | h_{\mathsf{d}} [m
  \chi_{\Omega}] |^2 . \label{eq:GLunorm}
\end{equation}
Here, $m \in H^1 (\Omega ; \Stwo^2)$, and $m \chi_{\Omega}$ \UUU is \EEE the extension by
zero of $m$ to $\R^3$. The {\tmem{Dirichlet energy}} $\mathcal{F}$, i.e., the
first term in {\eqref{eq:GLunorm}, penalizes nonuniformities in the
magnetization orientation, whereas the {\tmem{magnetostatic self-energy}}
$\mathcal{W}$, i.e., the second term in {\eqref{eq:GLunorm}, is the energy
associated with the demagnetizing field $h_{\text{d}}$ generated by $m
\chi_{\Omega}$, which describes the long-range dipole interaction of the
magnetic moments: for $v \in L^2 \left( \R^3, \R^3 \right)$, the demagnetizing
field $- h_{\mathsf{d}} [v] \in L^2 \left( \R^3, \R^3 \right)$ can be
characterized as the $L^2$-projection of $v$ on the space of gradients $\nabla
\dot{H}^1 \left( \R^3 \right) \assign \left\{ \nabla v \of v \in \mathcal{D}'
(\R^3), \nabla v \in L^2 (\R^3, \R^3) \right\}$,
see~{\cite{Di_Fratta_2019_var,praetorius2004analysis}} for details.

Much of the pattern observed in ferromagnetic materials is explained by the
competition between the two contributions in \eqref{eq:GLunorm}; in
particular, the formation of almost uniform magnetization regions
({\tmem{magnetic}} {\tmem{domains}}) separated by thin transition layers
({\tmem{domain walls}}), as predicted by the Weiss theory of
ferromagnetism~(cf.~{\cite{BrownB1963,hubert2008magnetic}}).

However, recent advancements in nanotechnology have led to the discovery of
{\tmem{magnetic skyrmions}}: chiral spin textures that carry a nontrivial
topological charge. Unlike conventional magnetic domains, magnetic skyrmions
exhibit unusual swirling textures, and they arise in ferromagnetic materials
with low crystallographic symmetry. The primary mechanism behind their
formation and stability is weak antisymmetric exchange interactions, also
known as {\tmem{Dzyaloshinskii--Moriya interactions}} (DMI), which result from
the combination of spin-orbit and superexchange
interactions~{\cite{dzyaloshinsky1958thermodynamic,moriya1960anisotropic}}.

In the continuum theory of micromagnetism, DMI is accounted through the
so-called {\tmem{chirality tensor}} $\nabla m \times m \assign (\partial_i m
\times m)_i$, whose components are the Lifshitz invariants of $m$ (see, e.g.,
{\cite[Supplementary information]{Hoffmann_2017}}). The {\tmem{bulk}} DMI
energy density corresponds to the trace of the chirality tensor: in the
presence of the bulk DMI, the micromagnetic energy functional
\eqref{eq:GLunorm} has to be remodeled under the form
\begin{equation}
  \mathcal{G} (m) \assign \mathcal{F} (m) +\mathcal{H}_{\tmop{bulk}} (m)
  +\mathcal{W} (m) \assign \frac{1}{2}  \int_{\Omega} | \nabla m|^2 + \gamma
  \int_{\Omega} \mathrm{curl} m \cdot m + \frac{1}{2}  \int_{\R^3} |
  h_{\mathsf{d}} [m \chi_{\Omega}] |^2 \underset{}{} . \label{eq:micromagDMI}
\end{equation}
The normalized constant $\gamma \in \R$ is the bulk DMI constant, and its sign
affects the chirality of the ferromagnetic
system~{\cite{bogdanov1994thermodynamically,nagaosa2013topological}} (see
{\cite[Sec.~4]{Di_Fratta_2023}} for further forms of DMI). Note that the bulk
DMI $\mathcal{H}_{\tmop{bulk}}$ in \eqref{eq:micromagDMI} is what the energy
term $\mathcal{H}$ in \eqref{func:loc0} reduces to when the Dzyaloshinskii
vectors have the form $\di \assign \gamma e_i$, with $e_i$ being the $i$-th
element of the standard basis of $\R^3$ (observe that $\tmop{curl} m = \sum_{i
= 1}^3 e_i \times \partial_i m$).

\subsection{The \texorpdfstring{$L^2$}{L2}-theory of symmetric and antisymmetric exchange
interactions. The Heisenberg setting}Ferromagnetism occurs in materials where
the spins tend to align with each other, thereby generating an observable
magnetic field outside of the media. {\tmem{Symmetric}} exchange interactions
are the primary mechanism behind this effect, and in the isotropic Heisenberg
model their behavior is described by the Hamiltonian
(cf.~{\cite{keffer1962moriya}})
\begin{equation}
  \mathcal{F}_{\Lambda} = - S^2  \sum_{(i, j) \in \Lambda \times \Lambda}
  J_{ij} m_i \cdot m_j . \label{eq:heisenbergmag}
\end{equation}
Here, $J_{ij} = J_{ji}$ is the symmetric exchange constant between the spins
$m_i, m_j \in \Stwo^2$, occupying the $i$-th and $j$-th site of the crystal
lattice $\Lambda$, and $S$ is the magnitude of the spin. The exchange constant
$J_{ij}$ is a positive quantity in ferromagnetism; it depends on the
ferromagnet's crystal structure and weighs the intensity of the interaction
among different spins.

When antisymmetric interactions cannot be neglected because, e.g., the
ferromagnetic crystal lacks inversion symmetry, the DMI induces a spin canting
of the magnetic moments, and the Hamiltonian \eqref{eq:heisenbergmag} has to
be remodeled under the form (cf.~{\cite{keffer1962moriya}})
\begin{equation}
  \mathcal{E}_{\Lambda} \assign \mathcal{F}_{\Lambda} +\mathcal{H}_{\Lambda}
  \assign - S^2  \sum_{(i, j) \in \Lambda \times \Lambda} J_{ij} m_i \cdot m_j
  + \sum_{(i, j) \in \Lambda \times \Lambda} d_{ij} \cdot (m_i \times m_j) .
  \label{eq:DMINaive}
\end{equation}
The {\tmem{Dzyaloshinskii vector}} $d_{ij} = - d_{ji}$ is an axial vector that
depends, other than from the relative distance between the spins $m_i$ and
$m_j$, on the symmetry of the crystal lattice; its precise form has to be
determined following Moriya's rules~{\cite{moriya1960anisotropic}}. We stress
that while the term $\mathcal{F}_{\Lambda}$ is symmetric, in the sense that
magnetic moments with right-handed $\left( { \left. \nwarrow
\hspace{0.17em} \nearrow \right)} \right.$ or left-handed $\left( {
\left. \nearrow \hspace{0.17em} \nwarrow \right)} \right.$ alignment give the
same contributions, the discrete antisymmetric exchange energy
$\mathcal{H}_{\Lambda}$ distinguishes between those two states via the local
chirality imposed by the Dzyaloshinskii vectors $d_{i  j}$.

In the limit of a continuous distribution of lattice sites occupying a region
$\Omega \subseteq \R^3$, the exchange energy \eqref{eq:DMINaive}, up to a
constant term, can be expressed as
\begin{align}
  \mathcal{E}_{J, d} (m) \assign \mathcal{F}_J (m) +\mathcal{H}_d (m) 
  & \assign
  \frac{1}{2}  \int_{\Omega \times \Omega} J (x - y)  |m (x) - m (y) |^2 \de x
  \, \de y \notag\\ 
  & \qquad\qquad\qquad + \int_{\Omega \times \Omega} d (x - y) \cdot (m (y) \times m (x)) 
  \de x \, \de y \label{eq:L2ExchDMI},
\end{align}
with $m : \Omega \rightarrow \Stwo^2$ the normalized magnetization density, $J
\geqslant 0$ a {\tmem{symmetric}} exchange kernel with even symmetry, i.e.,
such that $J (z) = J (- z)$ for every $z \in \R^3$, \ and $d$ an
{\tmem{antisymmetric}} exchange kernel  with odd symmetry, i.e., such that
$d (- z) = - d (z)$ for every $z \in \R^3$.

The energy functional $\mathcal{E}_{\varepsilon}$ in \eqref{func:nonloc} is a
faithful analog of the energy $\mathcal{E}_{J, d}$ in \eqref{eq:L2ExchDMI},
and we found it more convenient to work with $\mathcal{E}_{\varepsilon}$
because of scaling reasons. Specifically, we set $J (z) \assign \rho (z) / | z
|^2$ and $d (z) \assign \nu_{\varepsilon} (z) / | z |$. Overall, the exchange
terms $\mathcal{F}_{\ep}$, $\mathcal{H}_{\ep}$, in \eqref{eq:Fepssymexch} and
\eqref{eq:Hepsasymexch}, are the continuous counterparts, respectively, of the
symmetric Heisenberg Hamiltonian associated with a many-electron system and of
the antisymmetric exchange interactions due to the spin-orbit coupling between
neighboring magnetic spins.

\subsection{State of the art and contributions of the present work}At first
glance, the relationship between the nonlocal energies $\mathcal{F}_{\ep}$,
$\mathcal{H}_{\ep}$, in \eqref{eq:Fepssymexch}, \eqref{eq:Hepsasymexch}, and
the terms $\mathcal{F}$ and $\mathcal{H}$ in \eqref{func:loc0} is not evident,
but it can be formally revealed through a first-order asymptotic expansion of
the magnetization $m (x)$ in a neighborhood of $x \in \Omega$, i.e., by
setting $m (y) = m (x) + \mathD m (x) (y - x) +\mathcal{O} (| y - x |)$. The
asymptotic analysis \UUU becomes \EEE more reliable the more one can neglect variations
of $m$ around $x \in \Omega$, i.e., the more the kernels $\rho_{\varepsilon}$
and $\nu_{\varepsilon}$ concentrate their mass around the origin, i.e., the
more the exchange interactions act on a very short range only. The question
then is whether and in which sense the exchange energy $\mathcal{E}$ in
\eqref{func:nonloc} is a short-range approximation of the Heisenberg
$L^2$-description in \eqref{eq:L2ExchDMI}. In this paper, we give an
affirmative answer to \UUU this \EEE question under physically relevant hypotheses on the
kernels.

When only symmetric exchange interactions are considered, the (affirmative)
answer is already known because of the $\Gamma$-convergence result established
in {\cite{ponce2004new}}. Surprisingly, the motivation for the results in
{\cite{ponce2004new}} came from very different reasons: a question left open
in the seminal paper~{\cite{bourgain2001another}} where a new characterization
of Sobolev spaces is presented, and the Bourgain--Brezis--Mironescu (BBM)
formula made its first appearance. In~{\cite{bourgain2001another}},
{\tmem{pointwise}} convergence of $\mathcal{F}_{\ep}$ to $\mathcal{F}$ is
established when $(\rho_{\ep})_{\ep \in \R_+}$ is a family of radial
mollifiers concentrating their mass at the origin as $\ep \rightarrow 0$;
also, a technical lemma provides upper and lower bounds on the
{\tmem{variational}} convergence of $\mathcal{F}_{\ep}$, but not sufficient to
deduce $\Gamma$-convergence, finally obtained in {\cite{ponce2004new}}. In
summary, in the absence of DMI, the classical symmetric exchange energy can be
considered \UUU as \EEE the very short-rage limit of the family of nonlocal energies
$(\mathcal{F}_{\varepsilon})$, and the main aim of this paper is to extend the
result to the case in which also antisymmetric interactions are present and
under legitimate hypotheses on the antisymmetric exchange kernels
$\nu_{\varepsilon}$.

\textcolor{black}{Over the years, several papers have presented new BBM-type
formulas that expanded the original results in \cite{bourgain2001another} into various directions~\cite{brezis2016bbm,ponce2004new, MR2765717, MR3412379, MR1940355,MR2123016} (see
also~{\cite{squassina2016bourgain,pinamonti2017magnetic}} for some
applications to the magnetic Schr{\"o}dinger operator).} Also, formal
asymptotics predicts that the choice of appropriate exchange kernels is
crucial for $\mathcal{E}_{\varepsilon}$ to reduce to $\mathcal{E}$, \UUU and we
refer the reader to \cite{davoli2023sharp, foghem23} for some results about the
class of admissible kernels. \EEE

In contrast to the symmetric energy functional $\mathcal{F}_{\ep}$, the
variational convergence of $\mathcal{H}_{\ep}$ in the asymptotic regime of
very short-range interactions \UUU had \EEE not been investigated \UUU so far. \EEE To the best of the
authors' knowledge, no variational analysis involving $\mathcal{H}_{\ep}$ has
been \UUU carried out \EEE thus far; this includes basic questions like the existence
of $\Stwo^2$-valued minimizers of $\mathcal{E}_{\ep}$ or the qualitative
behaviors of minimizing sequences; these are significant issues because of
their direct bearing on the emergence of magnetic skyrmions, and \UUU will be the subject of a forthcoming paper. \EEE

\subsection*{\UUU Notation.} \UUU In what follows, we will adopt standard notation for the spaces of Lebesgue and Sobolev functions, as well as for the space of functions with bounded variation ($BV$). Balls with center $x$ and radius $r$ will be denoted by $B_r(x)$ and characteristic functions of sets $A$ by $\chi_A$, with the convention that $\chi_A(x)=1$ if $x\in A$ and $\chi_A(x)=0$ otherwise. Whenever not explicitly mentioned otherwise, $C>0$ will alway denote a generic constant, only dependent on the data of the problem, and whose value might, in principle, change from line to line. \EEE

\section{Statement of main results}\label{sec:rigsettingmainresults}
 For a given $\ep > 0$, we denote by $D_{\ep} (\Omega ; \Stwo^2)$ the subspace
of $L^2 (\Omega ; \Stwo^2)$ where $\mathcal{E}_{\ep}$ is finite:
\begin{equation}
  D_{\ep} (\Omega ; \Stwo^2) = \{m \in L^2 (\Omega ; \Stwo^2) :
  \mathcal{E}_{\ep} (m) < + \infty\} . \label{eq:defdomEeps}
\end{equation}
Also, since many of our results still hold when $m$ is not constrained to take
values on $\Stwo^2$, and even if $\Omega$ is unbounded, it is convenient for
us to denote by $D_{\ep} (\Omega ; \R^3)$ the analog unconstrained subspace of
$L^2 (\Omega ; \R^3)$. \textcolor{Black}{Of course, the space $D_{\ep} (\Omega ; \Stwo^2)$ on which the
energy $\mathcal{E}_{\ep}$ is finite depends on the particular choice of the
kernels $\rho_{\ep}$ and $\nu_{\ep}$. In what follows, we will make the
following assumptions.}

\subsection*{The symmetric exchange kernels \texorpdfstring{$\rho_{\varepsilon}$}{ρε}} Driven by
formal asymptotics, we consider a family of kernels $(\rho_{\ep})_{\ep}$
satisfying the following hypotheses already adopted in~{\cite{ponce2004new}}
---which are more general than the ones initially proposed in
{\cite{bourgain2001another}} (see also {\cite{davoli2023sharp}} for an
extensive discussion).

\begin{enumerate}[label=(G\arabic*)]
    \item \label{G1}For every $\ep > 0$, $\rho_{\ep} \geqslant 0$ in
  $\R^3$ and $\rho_{\ep} \in L^1 (\R^3)\UUU\cap BV_{\rm loc}(\R^3)\EEE$ with $\| \rho_{\ep} \|_{L^1 (\R^3)} =
  1$.  
  \item\label{G2}For every $\delta > 0$
  \begin{equation}
    \lim_{\ep \to 0}  \int_{\R^3 \setminus B_{\delta} (0)} \rho_{\ep} (y)
    \hspace{0.17em} \de y = 0.
  \end{equation}
  \item \label{G3}There exist linearly independent directions $v_1,
  v_2, v_3 \in \Stwo^2$ and $\delta > 0$ such that $C_{\delta} (v_i) \cap
  C_{\delta} (v_j) = \emptyset$ for $i \neq j$, and for any $i = 1, 2, 3$
  \begin{equation}
    \underset{\ep \to 0}{\limsup} \int_{C_{\delta} (v_i)} \rho_{\ep} (y) \de y
    > 0.
  \end{equation}
  Here, for $v \in \Stwo^2$ and $\delta > 0$, $C_{\delta} (v) \assign \R_+
  E_{\delta} (v)$ denotes the cone centered at the origin, whose projection on
  $\Stwo^2$ is given by $E_{\delta} (v) \assign \left\{ w \in \Stwo^2 \of w
  \cdot v > (1 - \delta) \right\}$.
\UUU\item \label{G4} There exists a constant $0<\kappa\leqslant 1$ such that, denoting by $\rho_\ep^{\rm rad}$ the radial functions $\rho_\ep^{\rm rad}:\R_+\to \R_+$, defined for every $\ep>0$ as
  \begin{equation}
  \label{eq:radial-k}
  \rho_\ep^{\rm rad}(t):={\rm ess}\inf\{\rho_\ep(x):\, x\in \R^3\text{ with }|x|=t\}\quad\text{ for every }t\in \R_+,
  \end{equation}
  there holds 
  \begin{equation}
  \label{eq:non-degenerate}
  \inf_{\ep>0}\int_{\R^3}\rho_\ep^{\rm rad}(|x|) \de x\geqslant \kappa.
  \end{equation}
  \end{enumerate}
\UUU

 \begin{remark}[On Hypotheses \ref{G1}, \ref{G3}, and \ref{G4}]
      \label{rk:G3-G4}
  A few words about Hypotheses \ref{G1},\ref{G3} and \ref{G4} are in order.
  \EEE
Hypothesis \ref{G3} is a weaker condition than the
\tmtextit{radiality} assumed in~{\cite{bourgain2001another}}. Roughly
speaking, it assures that in the limit of very short-range interactions ($\ep \to
0$), the family $(\rho_{\ep})_{\ep}$ has nontrivial support at least around
three linearly independent directions. Condition \ref{G3}
has been introduced in~{\cite{ponce2004new}} in order to prove a
characterization of $W^{1, p} (\Omega)$, $1 < p < \infty$. Our work uses it to
prove a regularity result for our nonlocal functional \eqref{func:nonloc} (see
Theorem \ref{thm:charact}). 
\UUU In view of \ref{G1}, Hypothesis \ref{G4} is automatically satisfied when the kernels $(\rho_\ep)_\ep$ are radial. If this is not the case, Hypothesis \ref{G4} provides a non-degeneracy condition for suitable \textcolor{Black}{radial lower envelopes} of the kernels. In particular, in the prototypical case in which the family $(\rho_\ep)_\ep$ is given by $$\rho_\ep(x):=\frac{1}{\ep^3}\rho\left(\frac{x}{\ep}\right),$$ for every $\ep>0$ and $x\in \R^3$, where  $\rho\in C^\infty_c(\R^3)$ is such that $0\leqslant \rho\leqslant 1$ in $\R^3$ and $\|\rho\|_{L^1(\R^3)}=1$, Hypothesis \ref{G4} is directly satisfied as long as \textcolor{Black}{$\rho$ is strictly positive at the origin}. In fact, assume that $\rho(0)=2\delta>0$. By the regularity of $\rho$, there exists a ball of radius $r>0$ such that $\rho(x)\geqslant \delta$ for every $x\in \overline{B_r(0)}$. In particular, $\rho_\ep^{\rm rad}(|x|)\geqslant \frac{\delta}{\ep^3}$ for every $\ep>0$ and \textcolor{Black}{$x\in \overline{B_{r\ep}(0)}$}, so that \eqref{eq:non-degenerate} holds with $\kappa=\delta|B_r(0)|$.

Finally, the regularity of the kernels in \ref{G1} can be weakened to just $L^1(\R^3)$ in the case in which the kernels are radial. The further $BV_{\rm loc}$- regularity is only needed to provide a meaning in the sense of traces to the restriction of the kernels $\rho_{\ep}$ to spheres of radius $r$, so that the definition in \eqref{eq:radial-k} is well-posed.  
\end{remark}
\EEE
 
\subsection*{The antisymmetric exchange kernels \texorpdfstring{$\nu_{\varepsilon}$}{vε}} Motivated
by formal asymptotics, we consider a family of vector-valued kernels
$(\nu_{\ep})_{\ep}$ satisfying the following assumptions.
\begin{enumerate}[label=(H\arabic*)]
  \item \label{HH1}For every $\ep > 0$, $\nu_{\ep}$ is odd, i.e.,
  $\nu_{\ep}  (- y) = - \nu_{\ep} (y)$ for each $y \in \R^3$, and $\nu_{\ep}
  \in L^1 (\R^3 ; \R^3)$ with $\| \nu_{\ep} \|_{L^1 (\R^3)} = 1$.
  
  \item \label{HH2}For every $\delta > 0$ there holds
  \begin{equation}
    \label{H2} \lim_{\ep \to 0}  \int_{\R^3 \setminus B_{\delta} (0)} |
    \nu_{\ep} (y) | \hspace{0.17em} \de y = 0.
  \end{equation}
  \item \label{HH3} For $i = 1, 2, 3$, the following limit exists
  \begin{equation}
    \lim_{\ep \to 0}  \int_{\R^3} \nu_{\ep} (y) \frac{y_i}{|y|}
    \hspace{0.17em} \de y = : \di \in \R^3 . \label{eq:wholeseqconvneps}
  \end{equation}
  Motivated by their physical significance, we will refer to the vectors $\di$
  as the Dzyaloshinskii vectors. Note that $| \di | \leqslant 1$ for every $i
  = 1, 2, 3$.
\end{enumerate}
Later on, we will need the following condition \ref{A1},
which connects the integrability of the two families of kernels
$(\rho_{\ep})_{\ep}$ and $(\nu_{\ep})_{\ep}$.
\begin{enumerate}[label=(A\arabic*)]
  \item \label{A1}There exist a constant $C > 0$ such that
  \begin{equation}
    \sup_{\ep > 0} \hspace{0.17em} \left\lVert \hspace{0.17em}
    \frac{\nu_{\ep}}{\rho_{\ep}} \hspace{0.17em} \right\rVert_{L^{\infty}
    (\R^3)} \leqslant C.
  \end{equation}
\end{enumerate}
Throughout the work, we specify in which cases assumption
\ref{A1} is needed.

\begin{remark}[\UUU On the odd symmetry of $\nu_\ep$\EEE]
  We emphasize that the requirement of $\nu_{\ep}$ being odd causes no loss of
  generality. Indeed, one can always decompose a vector-valued kernel into its
  even and odd parts, and an immediate symmetry argument shows that the
  presence of an even component in $\nu_{\ep}$ would give no contribution to
  $\mathcal{H}_{\ep}$ due to the antisymmetric nature of the functional
  $\mathcal{H}_{\ep}$ and its linearity with respect to the kernel. Also, our
  results still hold when one relaxes \ref{HH3}, requiring
  that \eqref{eq:wholeseqconvneps} holds up to the extraction of a
  subsequence, but we do not insist on these refinements.
\end{remark}

The main result \UUU of this paper \EEE consists of a compactness and $\Gamma$-convergence \UUU analysis \EEE
for our nonlocal exchange functionals in \eqref{func:nonloc}.

\begin{theorem}
  \textrm{{\opt Compactness and $\Gamma$-convergence\cpt \label{thm:gamma-conv}}} Let
  $\Omega \subseteq \R^3$ be a bounded Lipschitz domain. Assume
  \ref{G1}--\ref{G4},
  \ref{HH1}--\ref{HH3}, and
  \ref{A1}. The following assertions hold.
  \begin{enumerate}
    \item[(i)] {\tmname{(Compactness)}} If $(m_{\ep})_{\ep} \subset L^2 (\Omega ;
    \Stwo^2)$ is such that
    \begin{equation}
      \liminf_{\ep \rightarrow 0}  \hspace{0.17em} \mathcal{E}_{\ep} (m_{\ep})
      < + \infty, \label{eq:assumptionsupEeps}
    \end{equation}
    then, there exists $m \in H^1 (\Omega ; \Stwo^2)$ such that, possibly up
    to a non-relabeled subsequence, $m_{\varepsilon} \to m$ strongly in $L^2
    (\Omega ; \Stwo^2)$.
    
    \item[(ii)] {\tmname{($\Gamma$-convergence)}} There exists a finite Radon
    measure $\mu \in \mathcal{M} (\Stwo^2)$, with $\mu \left( \Stwo^2 \right)
    = 1$, such that, possibly up to a non-relabeled subsequence,
    \begin{equation}
      \Gamma \text{-} \lim_{\ep \to 0} \mathcal{E}_{\varepsilon}
      =\mathcal{E}_{\mu}, \quad \text{strongly in } L^2 \left( \Omega ;
      \Stwo^2 \right),
    \end{equation}
    with $\mathcal{E}_{\mu} (m) = + \infty$ if $m \in L^2 (\Omega ; \Stwo^2)
    \backslash H^1 (\Omega ; \Stwo^2)$ and
    \begin{equation}
      \mathcal{E}_{\mu} (m) = \int_{\Omega} \left( \int_{\Stwo^2} |
      \partial_{\sigma} m (x) |^2 \de \mu (\sigma) \right) \de x + \sum_{i =
      1}^3 \int_{\Omega} m (x) \cdot (d_i \times \partial_i m (x))
      \hspace{0.17em} \de x
    \end{equation}
    if $m \in H^1 (\Omega ; \Stwo^2)$.
  \end{enumerate}
\end{theorem}

\begin{remark}[\UUU The micromagnetic case\EEE]
  We note that choosing kernels $(\nu_{\ep})_{\ep}$ such that $\di = \gamma
  e_i$ for $i = 1, 2, 3$ in \ref{HH3}, and kernels $\left(
  \rho_{\ep} \right)_{\varepsilon}$ radial, from Theorem \ref{thm:gamma-conv}
  we obtain that with respect to the strong topology of $L^2$, there holds
  $\Gamma \text{-} \lim_{\ep \to 0} \mathcal{E}_{\ep} =\mathcal{E}$ with
  \begin{equation}
    \mathcal{E} (m) = \frac{1}{3} \int_{\Omega} | \nabla m (x) |^2  \de x +
    \gamma \int_{\Omega} \mathrm{curl} m (x) \cdot m (x) \hspace{0.17em} \de x
  \end{equation}
  if $m \in H^1 (\Omega ; \Stwo^2)$, and $\mathcal{E} (m) = + \infty$
  otherwise {\opt}if $m \in L^2 (\Omega ; \Stwo^2) \backslash H^1 (\Omega ;
  \Stwo^2)${\cpt}. Up to a constant factor, this is nothing but the
  micromagnetic energy functional \eqref{eq:micromagDMI} thoroughly studied in
  recent years as a theoretical foundation for the analysis of magnetic
  skyrmions emerging from bulk DMI.
\end{remark}

\UUU
\begin{remark}[Examples of kernels $\rho_\ep$ and $\nu_\ep$]
  The conditions imposed on the families of kernels $(\rho_{\ep})_{\ep}$ and
  $(\nu_{\ep})_{\ep}$ are not meant to be sharp but rather to cover most situations of interest in applications,
  in particular, for the analysis of magnetic skyrmions in chiral
  ferromagnetic materials. An explicit example of kernels satisfying the
  hypotheses above are presented below.  
  Let $\nu(y):=\frac{4}{|\Stwo^2|} y \chi_{B_1(0)}(y)$ and \textcolor{Black}{$\rho (y) := |\nu(y)|$}
  for every $y\in \R^3$, where we have denoted by $\chi_{B_1(0)}(y)$ the characteristic function of the unit ball in $\R^3$. Note that $\nu$ is odd, $\|\nu\|_{L^1(\R^3)}=\|\rho\|_{L^1(\R^3)}=1$, and $\rho$ is radial. Then, setting
  $$\nu_\ep(y):=\frac{1}{\ep^3}\nu\left(\frac{y}{\ep}\right)\quad\text{and }\rho_\ep(y):=\frac{1}{\ep^3}\rho\left(\frac{y}{\ep}\right)$$
  for every $\ep>0$ and $y\in \R^3$, Hypotheses \ref{G1}--\ref{G4}, as well as \ref{HH1}--\ref{HH2} and \ref{A1} are directly fulfilled. Additionally, \ref{HH3} holds with \textcolor{Black}{$d_i=\frac13 e_i$}, $i=1,2,3$. 
\end{remark}\EEE

\UUU
\begin{remark}[On the measure $\mu$]
Note that, on the one hand, as for the purely symmetric case, cf. {\cite[Lemma 8]{ponce2004new}}, the measure $\mu$ in Theorem \ref{thm:gamma-conv} (ii) is, in principle, dependent on the choice of the extracted subsequence. For some specific choices of the kernels $(\rho_\ep)_\ep$, on the other hand, for example when they are radial, the limiting functional is uniquely identified independently of the extracted subsequence, see also Remark \ref{rk:quadratic} below, so that the Gamma-convergence result actually holds for the entire sequence $(\mathcal{E}_{\ep})_{\ep}$.
\end{remark} \EEE

The proof of Theorem~\ref{thm:gamma-conv} is \UUU provided \EEE in
Section~\ref{sec:Gammaconv}. The argument \UUU relies upon \EEE the following results, which
are of interest in their own rights. The first assures that the nonlocal
functional \eqref{func:nonloc} is well-defined on $H^1 (\Omega ; \Stwo^2)$.

\begin{theorem}
  \label{thm:charact}Assume \ref{G1} and
  \ref{HH1}. For every $m \in H^1 (\Omega ; \R^3)$ there
  holds
  \begin{equation}
    \sup_{\varepsilon > 0} \mathcal{E}_{\ep} (m) \leqslant C \hspace{0.17em}
    \|m\|^2_{H^1 (\Omega)}, \label{tobeprovedreg1}
  \end{equation}
  for some constant $C > 0$ depending only on $\Omega$.
  
  Viceversa, assume \ref{G1}--\ref{G3},
  as well as \ref{HH1} and \ref{A1}. If
  $m \in L^2 (\Omega ; \R^3)$ is such that
  \begin{equation}
    \sup_{\ep > 0} \hspace{0.17em} \mathcal{E}_{\ep} (m) < + \infty,
    \label{eq:sup}
  \end{equation}
  then $m \in H^1 (\Omega ; \R^3)$.
\end{theorem}

\begin{remark}[\UUU The domain of the nonlocal functionals\EEE]
  Estimate \eqref{tobeprovedreg1} gives, in particular
  ({\tmabbr{cf.}}~\eqref{eq:defdomEeps}), that $H^1 (\Omega ; \R^3) \subseteq
  D_{\ep} (\Omega ; \R^3)$ and $H^1 (\Omega ; \Stwo^2) \subseteq D_{\ep}
  (\Omega ; \Stwo^2)$. Also, it will be evident from the proof that $H^1
  (\Omega ; \R^3) \subseteq D_{\ep} (\Omega ; \R^3)$ even when $\Omega$ is
  unbounded. In particular, $H^1 (\R^3 ; \R^3) \subseteq D_{\ep} (\R^3 ;
  \R^3)$.
\end{remark}

\begin{remark}[\UUU On the assumptions of Theorem \ref{thm:charact}\EEE]
  Estimate \eqref{tobeprovedreg1} guarantees that if
  \ref{G1} and \ref{HH1} hold and $m \in
  H^1 (\Omega ; \R^3)$, then $\sup_{\ep > 0} \mathcal{E}_{\ep} (m) < +
  \infty$. The converse statement holds under the additional assumptions
  \ref{G2}, \ref{G3}, 
  \ref{HH1}, and
  \ref{A1}. Assumptions \ref{G2}--\ref{G3} are the main \UUU requirements \EEE in
  {\cite[Theorem~5]{ponce2004new}}. Roughly speaking,
  \ref{G3} prevents \EEE degenerate mass behavior of the family
  $(\rho_{\ep})_{\ep}$, so that the partial derivatives of $m$ can be
  controlled over all directions. The further assumption
  \ref{A1} expresses an integrability relation between the
  two families of kernels $(\rho_{\ep})_{\ep}$ and $(\nu_{\ep})_{\ep}$ which,
  nevertheless, naturally arise when modeling antisymmetric exchange
  interactions in the variational theory of micromagnetism.
\end{remark}

Another crucial ingredient \UUU for \EEE the proof of Theorem~\ref{thm:gamma-conv}
is a uniform convergence result for the antisymmetric exchange energies
$\mathcal{H}_{\ep}$, contained in the following Theorem~\ref{thm:point-conv}.

\begin{theorem}
  {\tmname{\opt Pointwise and uniform convergence of the antisymmetric exchange
  term\cpt}}\label{thm:point-conv} Assume
  \ref{HH1}--{\hyperrefstar{(H3)}{H3}} and \UUU let $\Omega$ be \EEE
  bounded. For every $m \in H^1 (\Omega ; \R^3)$ there holds
  \begin{equation}
    \lim_{\ep \to 0} \mathcal{H}_{\ep} (m) =\mathcal{H} (m),
    \label{eq:Hepspointwiseconv}
  \end{equation}
  where, we recall {\tmem{(}}see {\tmem{\eqref{eq:Fepssymexch}}} and
  {\tmem{\eqref{func:loc0})}}
  \[ \mathcal{H}_{\ep} (m) = \iint_{\Omega \times \Omega} \frac{\nu_{\ep}  (x
     - y)}{|x - y|} \cdot (m (x) \times m (y)) \hspace{0.17em} \de x
     \hspace{0.17em} \de y \hspace{0.17em}, \quad \mathcal{H} (m) = \sum_{i =
     1}^3 \int_{\Omega} m (x) \cdot (d_i \times \partial_i m (x))
     \hspace{0.17em} \de x, \]
  and the $d_i \in \R^3$ are the Dzyaloshinskii vectors defined in
  {\tmem{\ref{HH3}}}.
  
  Moreover, if $(m_{\varepsilon})_{\varepsilon > 0}$ is a family in $C^2 
  (\bar{\Omega} ; \R^3)$ such that $m_{\varepsilon} \rightarrow m$ in $C^2 
  (\bar{\Omega} ; \R^3)$, then
  \begin{equation}
    \lim_{\ep \to 0} \mathcal{H}_{\ep} (m_{\varepsilon}) =\mathcal{H} (m) .
    \label{eq:Hepspointwiseconvunif}
  \end{equation}
\end{theorem}

\begin{remark}[\UUU The case of unbounded $\Omega$\EEE]
  It will be evident from the proof that \eqref{eq:Hepspointwiseconv} holds
  even when $\Omega$ is unbounded. In particular, when $\Omega = \R^3$. Also,
  for simplicity, we required \ref{HH3} to ensure the
  convergence of the whole sequence in \eqref{eq:wholeseqconvneps}. It will be
  clear from the proof that this requirement can be weakened by only requiring
  that the convergence in \ref{HH3} holds up to the
  extraction of a subsequence.
\end{remark}

As a direct consequence of Theorem \ref{thm:point-conv}, and {\cite[Theorem
1]{ponce2004new}}, we \UUU infer \EEE the following result about the pointwise asymptotic
behavior of the total nonlocal energy \eqref{func:nonloc}.

\begin{corollary}
  {\tmname{\opt Pointwise convergence of the total energy\cpt}}\label{cor:point-conv}
  Assume
  {\tmem{\ref{G1}}}--{\tmem{\ref{G3}}},
  as well as
  \ref{HH1}--{\tmem{\ref{HH3}}} \UUU and \ref{A1}. Then, there \EEE
  exists a finite Radon measure $\mu \in \mathcal{M} (\Stwo^2)$, with $\mu
  \left( \Stwo^2 \right) = 1$, such that, possibly up to a non-relabeled
  subsequence, for every $m \in H^1 (\Omega ; \R^3)$, there holds
  \begin{equation}
    \lim_{\varepsilon \rightarrow 0} \mathcal{E}_{\varepsilon} (m) =
    \int_{\Omega} \left( \int_{\Stwo^2} | \partial_{\sigma} m|^2  \de \mu
    (\sigma) \right)  \de x + \sum_{i = 1}^3 \int_{\Omega} m (x) \cdot (d_i
    \times \partial_i m (x)) \hspace{0.17em} \de x, \label{eq:Poncelimit}
  \end{equation}
  with the $d_i \in \R^3$ being the Dzyaloshinskii vectors introduced in
  {\tmem{\ref{HH3}}}.
  
  Moreover, if the family $(\rho_{\varepsilon})_{\varepsilon}$ consists of
  radial kernels, then for every $m \in H^1 (\Omega ; \R^3)$ we have
  \[ \lim_{\varepsilon \rightarrow 0} \left[ \mathcal{F}_{\varepsilon} (m)
     +\mathcal{H}_{\ep} (m) \right] = \frac{1}{3} \int_{\Omega} | \nabla m|^2
     \hspace{0.17em} \de x + \sum_{i = 1}^3 \int_{\Omega} m (x) \cdot (d_i
     \times \partial_i m (x)) \hspace{0.17em} \de x. \]
\end{corollary}

\begin{remark}[\UUU Structure of the limiting symmetric term\EEE]\label{rk:quadratic}
  Since we are in the quadratic setting, one has $| \partial_{\sigma} m|^2 =
  \sum_k (\sigma \otimes \sigma) \nabla m_k \cdot \nabla m_k$ and, therefore,
  one can rewrite the first term in \eqref{eq:Poncelimit} under the form
  \begin{equation}
    \lim_{\varepsilon \rightarrow 0} \mathcal{F}_{\varepsilon} (m) = \sum_k
    \int_{\Omega} A \nabla m_k \cdot \nabla m_k  \de x
  \end{equation}
  where $A$ is the anisotropic matrix given by $A \assign \int_{\Stwo^2}
  (\sigma \otimes \sigma) \de \mu (\sigma)$. If the kernels are radial, then
  $\mu$ is isotropic, in the sense that the resulting matrix $A$ is given by
  $(1 / 3) I$ with $I$ being the $3 \times 3$ identity matrix.
\end{remark}

\section{Pointwise and uniform behavior of the energy (proofs of
Theorem~\ref{thm:charact} and Theorem~\ref{thm:point-conv})}\label{Sec:Point}

\begin{proof}[\tmname{\bf Proof of {Theorem~\ref{thm:charact}}}]
  We split the proof in two steps. In {\tmname{Step~1}}, we derive the
  estimate \eqref{tobeprovedreg1}, while in {\tmname{Step~2}}, we prove that
  if \eqref{eq:sup} holds, then $m \in H^1 (\Omega ; \R^3)$.{\medskip}
  
  {\noindent}{\tmname{Step~1}}. Recall that, by definition, $\mathcal{E}_{\ep}
  (m) =\mathcal{F}_{\ep} (m) +\mathcal{H}_{\ep} (m)$. For ease of computation,
  we treat the two energy terms separately. For $m \in H^1 (\Omega ; \R^3)$,
  we denote by $\tilde{m} \in H^1 (\R^3 ; \R^3)$ the extension of $m$ to the
  whole of $\R^3$ given by the classical extension operator on $H^1 (\Omega ;
  \R^3)$.
  
  For the symmetric exchange term $\mathcal{F}_{\ep}$ we have
  
  \begin{align}
    \nonumber \mathcal{F}_{\ep} (m) & \leqslant \iint_{\R^3 \times \R^3} \rho_{\ep}  (x
    - y)  \frac{| \tilde{m} (x) - \tilde{m} (y) |^2}{|x - y|^2}
    \hspace{0.17em} \de x \hspace{0.17em} \de y \\
    & = \int_{\R^3} \rho_{\ep} (h) \left( \int_{\R^3} \frac{| \tilde{m} (x -
    h) - \tilde{m} (x) |^2}{|h|^2} \hspace{0.17em} \de x \right) \de h 
    \label{change}\\
    \nonumber & = \int_{\R^3} \rho_{\ep} (h) \hspace{0.17em} \frac{\| \tau_{- h} 
    \tilde{m} - \tilde{m} \|^2_{L^2 (\R^3)}}{|h|^2} \hspace{0.17em} \de h \\
    \nonumber & \leqslant \| \rho_{\ep} \|_{L^1 (\R^3)}  \hspace{0.17em} \| \nabla
    \tilde{m} \|^2_{L^2 (\R^3)} \\
     & \leqslant \hspace{0.17em} C \hspace{0.17em} \|m\|^2_{H^1 (\Omega)}, 
    \label{Brezis}
  \end{align}
  
  for some constant $C > 0$ depending only on $\Omega$ through the linear
  extension operator. The equality \eqref{change} is the result of the change
  of variables $y \mapsto x - h$ for fixed $x \in \R^3$, whereas to get
  \eqref{Brezis}, we used \ref{G1} and standard properties
  of Sobolev spaces (cf., e.g.,~{\cite[Thm.~1]{bourgain2001another}} or
  {\cite[Prop.~9.3]{brezis2011functional}}).
  
  Similarly, for the antisymmetric exchange term, we obtain
  
  \begin{align}
     |\mathcal{H}_{\ep} (m) | & \leqslant \iint_{\Omega \times
    \Omega} \frac{\left| \nu_{\ep} (x - y) \right|}{|x - y|} \cdot | m (x)
    \times [m (y) - m (x)] | \hspace{0.17em} \de x \hspace{0.17em} \de y
    \hspace{0.17em} \nonumber\\
    \nonumber & \leqslant \iint_{\R^3 \times \R^3} \frac{| \nu_{\ep} (x - y) |}{|x -
    y|}  | \tilde{m} (x) || \tilde{m} (y) - \tilde{m} (x) | \hspace{0.17em}
    \de x \hspace{0.17em} \de y \\
    & = \int_{\R^3} | \nu_{\ep} (h) | \left( \int_{\R^3} | \tilde{m} (x) |
    \frac{| \tilde{m} (x + h) - \tilde{m} (x) |}{|h|} \hspace{0.17em} \de x
    \right) \de h,  \label{h}
  \end{align}
  
  where in \eqref{h} we performed the change of variables $y \mapsto x + h$ for
  fixed $x \in \R^3$ and then used the odd symmetry of the kernel $\nu_{\ep}$.
  Applying H{\"o}lder's inequality as well as classical properties of Sobolev
  spaces (see~{\cite[Prop.~9.3]{brezis2011functional}}), from \eqref{h} and
  \ref{HH1} we infer that
  
  \begin{align*}
    |\mathcal{H}_{\ep} (m) | & \leqslant \int_{\R^3} | \nu_{\ep} (h) | \|
    \tilde{m} \|_{L^2 (\R^3)} \frac{\| \tau_h  \tilde{m} - \tilde{m} \|_{L^2
    (\R^3)}}{|h|} \hspace{0.17em} \de h \\
    & \leqslant \| \nu_{\ep} \|_{L^1 (\R^3)} \| \tilde{m} \|_{L^2 (\R^3)}  \|
    \nabla \tilde{m} \|_{L^2 (\R^3)} \\
    & \leqslant C \hspace{0.17em} \|m\|^2_{H^1 (\Omega)}, 
  \end{align*}
  
  for some constant $C > 0$ depending only on $\Omega$. This completes the
  proof of \eqref{tobeprovedreg1}.{\medskip}
  
  {\noindent}{\tmname{Step~2}}. We claim that it is sufficient to show that if
  $\sup_{\ep > 0} \hspace{0.17em} \mathcal{E}_{\ep} (m) < + \infty$, then
  $\sup_{\ep > 0} \hspace{0.17em} \mathcal{F}_{\ep} (m) < + \infty$. Indeed,
  as soon as we show that $\sup_{\ep > 0} \hspace{0.17em} \mathcal{F}_{\ep}
  (m) < + \infty$, we can then invoke {\cite[Theorem~5]{ponce2004new}} which
  shows that if $\sup_{\ep > 0} \hspace{0.17em} \mathcal{F}_{\ep} (m) < +
  \infty$ and the family $\left( \rho_{\ep} \right)_{\varepsilon}$ satisfies
  \ref{G1}--\ref{G3}, then $m \in H^1
  (\Omega ; \R^3)$ and there exists $\alpha > 0$ (depending only on
  the family of kernels $\left( \rho_{\ep} \right)_{\varepsilon}$) such that
  \begin{equation}
    \alpha^2 \int_{\Omega} | \nabla m|^2  \de x \leqslant \limsup_{\ep \to 0}
    \iint_{\Omega \times \Omega} \rho_{\ep}  (x - y)  \frac{|m (x) - m (y)
    |^2}{|x - y|^2} \hspace{0.17em} \de x \hspace{0.17em} \de y.
    \label{ineq:alpha}
  \end{equation}
  Since $\mathcal{H}_{\ep} (m)$ is a quantity without a definite sign, we use
  Young's inequality for products to estimate the term $\mathcal{H}_{\ep} (m)$
  from below. We get that for every $\delta > 0$ there holds
  \begin{align}
    \mathcal{H}_{\ep} (m) & =  \iint_{\Omega \times \Omega} \frac{m (x) - m
    (y)}{|x - y|} \cdot \left( m (y) \times \nu_{\ep} (x - y) \right)
    \hspace{0.17em} \de x \hspace{0.17em} \de y \nonumber\\
    & =  \iint_{\Omega \times \Omega} \left( \sqrt{\rho_{\ep} (x - y)}
    \frac{m (x) - m (y)}{|x - y|} \right) \cdot \left( m (y) \times
    \frac{\nu_{\ep} (x - y)}{\sqrt{\rho_{\ep} (x - y)}} \right)
    \hspace{0.17em} \de x \hspace{0.17em} \de y \nonumber\\
    & \geqslant  - \delta^2 \mathcal{F}_{\ep} (m) - \frac{1}{4 \delta^2}
    \iint_{\Omega \times \Omega}  \frac{|m (y) \times \nu_{\ep} (x - y)
    |^2}{\rho_{\ep}  (x - y)} \de x \hspace{0.17em} \de y. 
  \end{align}
  Adding $\mathcal{F}_{\ep} (m)$ to both sides of the previous estimate, we
  get that for every $\delta > 0$ there holds
  \begin{equation}
    (1 - \delta^2) \mathcal{F}_{\ep} (m) - \frac{1}{4 \delta^2} \iint_{\Omega
    \times \Omega}  \frac{|m (y) \times \nu_{\ep} (x - y) |^2}{\rho_{\ep}  (x
    - y)} \de x \hspace{0.17em} \de y  \leqslant  \mathcal{E}_{\ep} (m) . 
    \label{eq:square}
  \end{equation}
  Now, by \ref{A1}, and taking into account that $\|
  \nu_{\ep} \|_{L^1 (\R^3)} = 1$, we infer
  \begin{align}
    \iint_{\Omega \times \Omega}  \frac{|m (y) \times \nu_{\ep} (x - y)
    |^2}{\rho_{\ep}  (x - y)} \de x \hspace{0.17em} \de y & \leqslant 
    \int_{\Omega} |m (y) |^2 \left( \int_{\R^3} \frac{| \nu_{\ep} (z)
    |^2}{\rho_{\ep} (z)} \de z \right) \de y \nonumber\\
    & \leqslant  \frac{1}{4} \left( \int_{\Omega} |m (y) |^2 \de y \right)
    \lVert \nu_{\ep} \rVert_{L^1 (\R^3)} \hspace{0.17em} \left\lVert
    \frac{\nu_{\ep}}{\rho_{\ep}} \right\rVert_{L^{\infty} (\R^3)} \nonumber\\
    & \leqslant  C \hspace{0.17em} \|m\|_{L^2 (\Omega)}^2,
    \label{dis:L2control}
  \end{align}
  with $C > 0$ proportional to the constant from \ref{A1}.
  Hence, from \eqref{eq:square} and \eqref{dis:L2control} we deduce
  \begin{equation}
    (1 - \delta^2) \mathcal{F}_{\ep} (m) \leqslant \mathcal{E}_{\ep} (m) + C
    \|m\|_{L^2 (\Omega)}^2 \label{eq:funestfepshepsineeps}
  \end{equation}
  By taking $\delta^2 \assign 1 / 2$ we conclude that if $\sup_{\ep > 0}
  \hspace{0.17em} \mathcal{E}_{\ep} (m) < + \infty$ then $\sup_{\ep > 0}
  \hspace{0.17em} \mathcal{F}_{\ep} (m) < + \infty$ and hence, by
  {\cite[Theorem~5]{ponce2004new}}, that $m \in H^1 (\Omega ; \R^3)$.
\end{proof}

\begin{proof}[\tmname{\bf Proof of Theorem~\ref{thm:point-conv}}]
  \label{sub:point-conv}We split the proof into two main steps. In
  {\tmname{Step~1}}, we prove \eqref{eq:Hepspointwiseconvunif}, which, in
  particular, gives us \eqref{eq:Hepspointwiseconv} under the additional
  hypothesis that $m \in C^2  (\bar{\Omega} ; \R^3)$. Then, in
  {\tmname{Step~2}}, we use a density argument to show that if
  \eqref{eq:Hepspointwiseconv} holds for every $m \in C^2  (\bar{\Omega} ;
  \R^3)$ then it holds for every $m \in H^1 (\Omega ; \R^3)$.{\medskip}
  
  {\noindent}{\tmname{Step~1}}. Let $y \in \Omega$. For $R > 0$ such that $R
  \leqslant \text{dist} (y, \partial \Omega)$, we split $\mathcal{H}_{\ep}
  (m_{\varepsilon})$ as
  \begin{align}
    \mathcal{H}_{\ep} (m_{\varepsilon}) & =  \int_{\Omega} \left( \int_{B_R
    (y)} \frac{\nu_{\ep} (x - y)}{|x - y|} \cdot (m_{\varepsilon} (x) \times
    m_{\varepsilon} (y)) \hspace{0.17em} \de x \right)  \de y \nonumber\\
    &   \qquad + \int_{\Omega} \left( \int_{\Omega \setminus B_R (y)}
    \frac{\nu_{\ep} (x - y)}{|x - y|} \cdot (m_{\varepsilon} (x) \times
    m_{\varepsilon} (y)) \hspace{0.17em} \de x \right)  \de y \notag \\
    & \backassign  \int_{\Omega} I_{\ep, R} (m_{\varepsilon}, y) + J_{\ep,
    R} (m_{\varepsilon}, y)  \de y  \label{split:ballunif}
  \end{align}
  Outside of the ball $B_R (y)$, we note that
  \begin{align*}
    \left| J_{\ep, R} (m_{\varepsilon}, y) \right|  \leqslant & \int_{\Omega
    \setminus B_R (y)} \frac{| \nu_{\ep} (x - y) |}{|x - y|}  |m_{\varepsilon}
    (x) ||m_{\varepsilon} (y) | \hspace{0.17em} \de x \\
    & \leqslant  \frac{1}{R} \|m_{\varepsilon} \|^2_{L^{\infty} (\Omega)} 
    \int_{\R^3 \setminus B_R (0)} | \nu_{\ep} (x) | \hspace{0.17em} \de x. 
  \end{align*}
  Sending $\ep \to 0$ (for fixed $R > 0$), by assumption
  \ref{HH2} we find that $\lim_{\ep \to 0}  \|
  \nu_{\varepsilon} \|_{L^1 \left( \R^3 \setminus B_R (0) \right)} = 0$ and,
  therefore, that
  \begin{equation}
    \lim_{\ep \to 0} \, \left| J_{\ep, R} (m_{\varepsilon}, y) \right| = 0.
    \label{eq:B_omunif}
  \end{equation}
  We now take into account the contribution to $\mathcal{H}_{\ep}
  (m_{\varepsilon})$ inside of the ball $B_R (y)$, i.e., the term $I_{\ep}
  (m_{\varepsilon}, y)$ in \eqref{split:ballunif}. Given that
  $m_{\varepsilon} \in C^2  (\bar{\Omega} ; \R^3)$, we perform a
  first-order Taylor expansion inside $B_R (y)$ of the form $m_{\varepsilon}
  (x) = m_{\varepsilon} (y) + \mathD m_{\varepsilon} (y) (x - y) +
  R_{\varepsilon} (x, y)$ with (up to a constant factor)
  \[ | R_{\varepsilon} (x, y) | \leqslant \| D^{(2)} m_{\varepsilon}
     \|_{L^{\infty}} \cdot | x - y |^2 \]
  to find
  \begin{align}
    I_{\ep, R} (m_{\varepsilon}, y) & =  \int_{B_R (y)} \frac{\nu_{\ep}  (x -
    y)}{|x - y|} \cdot (\mathD m_{\varepsilon} (y) (x - y) \times
    m_{\varepsilon} (y)) \hspace{0.17em} \de x  \label{firstunif}\\
    \nonumber &   \qquad \qquad + \int_{B_R (y)} \frac{\nu_{\ep}  (x - y)}{|x - y|}
    \cdot (R_{\varepsilon} (x, y) \times m_{\varepsilon} (y)) \hspace{0.17em}
    \de x  \\
    \nonumber & \backassign  I_{\ep, R}^{(1)} (m_{\varepsilon}, y) + I_{\ep, R}^{(2)}
    (m_{\varepsilon}, y) . 
  \end{align}
  Let us focus on the \UUU first \EEE term on the right-hand side of \eqref{firstunif}. Simple
  algebra gives that
  \begin{align}
    I_{\ep, R}^{(1)} (m_{\varepsilon}, y) & =  m_{\varepsilon} (y) \cdot
    \int_{B_R (y)} \left( \frac{\nu_{\ep} (x - y)}{|x - y|} \times \partial_i
    m_{\varepsilon} (y) (x_i - y_i) \right) \hspace{0.17em} \de x \nonumber\\
    & =  \sum_{i = 1}^3 \left( \int_{B_R (y)} \frac{(x_i - y_i)}{|x - y|}
    \nu_{\ep} (x - y) \hspace{0.17em} \de x \right) \cdot (\partial_i
    m_{\varepsilon} (y) \times m_{\varepsilon} (y)) . 
    \label{eq:tobesplitrhs1unif}
  \end{align}
  Next, we observe that as a consequence of \ref{HH2} we
  have
  \begin{equation}
    \lim_{\ep \to 0} \left| \int_{\R^3 \setminus B_R (y)} \frac{\nu_{\ep} (x -
    y)}{|x - y|} (x_i - y_i) \hspace{0.17em} \de x \hspace{0.17em} \right|
    \leqslant \lim_{\ep \to 0}  \int_{\R^3 \setminus B_R (y)} | \nu_{\ep} (x -
    y) | \hspace{0.17em} \de x = 0. \label{eq:thanksH2unif}
  \end{equation}
  Therefore, combining \eqref{eq:thanksH2unif} and \ref{HH3}
  we infer that for $i = 1, 2, 3$, there holds
  \begin{align}
    \lim_{\ep \to 0}  \int_{B_R (y)} \frac{(x_i - y_i)}{|x - y|} \nu_{\ep} (x
    - y) \hspace{0.17em} \de x & =  \lim_{\ep \to 0}  \int_{\R^3} \frac{(x_i
    - y_i)}{|x - y|} \nu_{\ep} (x - y) \hspace{0.17em} \de x \nonumber\\
    & =  \lim_{\ep \to 0}  \int_{\R^3} \frac{x_i}{|x|} \nu_{\ep} (x)
    \hspace{0.17em} \de x \nonumber\\
    & =  \di .  \label{eq:diunif}
  \end{align}
  As a consequence, from \eqref{eq:tobesplitrhs1unif} and \eqref{eq:diunif} we
  conclude 
  \begin{equation}
    \lim_{\ep \to 0} I_{\ep, R}^{(1)} (m_{\varepsilon}, y) = \sum_{i = 1}^3 m
    (y) \cdot (d_i \times \partial_i m (y)) \label{eq:A1_omunif} .
  \end{equation}
  Note that \eqref{eq:A1_omunif} does not depend on $R$.
  
  For the term $I_{\ep, R}^2 (m_{\varepsilon}, y)$, we
  observe that
  \begin{align}
    \left| I_{\ep, R}^{(2)} (m_{\varepsilon}, y) \right| & \leqslant  |
    m_{\varepsilon} (y) | \int_{B_R (y)} | \nu_{\ep} (x - y) |
    \frac{|R_{\varepsilon} (x, y) |}{|x - y|} \hspace{0.17em} \de x
    \nonumber\\
    & \leqslant  R | m_{\varepsilon} (y) |  \| D^{(2)} m_{\varepsilon}
    \|_{L^{\infty}} \cdot \int_{\R^3} | \nu_{\ep} (x) | \hspace{0.17em} \de x
    \nonumber\\
    & =  R | m_{\varepsilon} (y) |  \| D^{(2)} m_{\varepsilon}
    \|_{L^{\infty}} . \nonumber
  \end{align}
  For the latter equality we used assumption \ref{HH1}.
  Thus, \UUU since $m_\ep\to m$ strongly in $C^2(\overline \Omega; \R^3)$, \EEE
  \begin{equation}
    \lim_{R \rightarrow 0} \lim_{\ep \to 0} \hspace{0.17em} \left| I_{\ep,
    R}^{(2)} (m_{\varepsilon}, y) \right| = 0. \label{eq:limRto0unif}
  \end{equation}
  Overall, combining estimates \eqref{eq:B_omunif},
  \eqref{eq:tobesplitrhs1unif}, \eqref{eq:diunif}, \eqref{eq:A1_omunif}, and
  \eqref{eq:limRto0unif}, we \UUU infer \EEE
  \begin{align}
    \lim_{\varepsilon \rightarrow 0} \int_{\Omega} \frac{\nu_{\ep}  (x -
    y)}{|x - y|} \cdot (m_{\varepsilon} (x) \times m_{\varepsilon} (y))
    \hspace{0.17em} \de x & =  \lim_{R \rightarrow 0} \lim_{\varepsilon
    \rightarrow 0}  \left( I_{\ep, R} (m_{\varepsilon}, y) + J_{\ep, R}
    (m_{\varepsilon}, y) \right)  \label{eq:isdomconvokunif}\\
    \nonumber & =  \lim_{R \rightarrow 0} \lim_{\varepsilon \rightarrow 0} I_{\ep, R}
    (m_{\varepsilon}, y) \\
    \nonumber & =  \lim_{R \rightarrow 0} \lim_{\varepsilon \rightarrow 0}  \left(
    I_{\ep, R}^{(1)} (m_{\varepsilon}, y) + I_{\ep, R}^{(2)} (m_{\varepsilon},
    y) \right) \\
    \nonumber & =  \lim_{\varepsilon \rightarrow 0} I_{\ep, R}^{(1)} (m_{\varepsilon},
    y) \\
    & =  \sum_{i = 1}^3 m (y) \cdot (d_i \times \partial_i m (y)) . 
    \label{eq:essentialconvHepsunif}
  \end{align}
  \UUU In view of \EEE \eqref{eq:essentialconvHepsunif}, to prove
  \begin{equation}
    \lim_{\varepsilon \rightarrow 0} \mathcal{H}_{\ep} (m_{\varepsilon}) =
    \sum_{i = 1}^3 \int_{\Omega} m (y) \cdot (d_i \times \partial_i m (y)) 
    \de y,
  \end{equation}
  it is sufficient to show that the family of functions on the left-hand side
  of \eqref{eq:isdomconvokunif} falls under the hypotheses of the Dominated
  Convergence Theorem. But this is indeed the case because of \UUU the estimate \EEE
  \begin{align}
    \left| \int_{\Omega} \frac{\nu_{\ep} (x - y)}{|x - y|} \cdot
    (m_{\varepsilon} (x) \times m_{\varepsilon} (y)) \hspace{0.17em} \de x
    \right| & \leqslant  \|m_{\varepsilon} \|_{L^{\infty} (\Omega)} 
    \int_{\Omega} | \nu_{\ep} (x - y) | \frac{|m_{\varepsilon} (x) -
    m_{\varepsilon} (y) |}{|x - y|} \hspace{0.17em} \de x \nonumber\\
    & \leqslant  L_{m_{\varepsilon}} \hspace{0.17em} \|m_{\varepsilon}
    \|_{L^{\infty} (\Omega)} \int_{\R^3} | \nu_{\ep} (x) | \hspace{0.17em} \de
    x \nonumber\\
    & =  L_{m_{\varepsilon}} \hspace{0.17em} \|m_{\varepsilon}
    \|_{L^{\infty} (\Omega)}, 
  \end{align}
  where $L_{m_{\varepsilon}}$ is the Lipschitz constant of $m_{\varepsilon}$,
  which is uniformly bounded in $\varepsilon$.{\medskip}
  
  {\noindent}{\tmname{Step~2}}. This second step concludes the proof by
  showing, through a density argument, that if \eqref{eq:Hepspointwiseconv}
  holds for every $m \in C^2  (\bar{\Omega} ; \R^3)$ then it holds for every
  $m \in H^1 (\Omega ; \R^3)$.
  
  We note that, for every $m_1, m_2 \in H^1 (\Omega ; \R^3)$, there holds
  \begin{align}
    \mathcal{H}_{\ep} (m_1) -\mathcal{H}_{\ep} (m_2) & =  \iint_{\Omega
    \times \Omega} \frac{\nu_{\ep}  (x - y)}{|x - y|} \cdot ((m_1 (x) - m_1
    (y)) \times m_1 (y)) \de x \hspace{0.17em} \de y \nonumber\\
    &   \qquad - \iint_{\Omega \times \Omega} \frac{\nu_{\ep}  (x - y)}{|x -
    y|} \cdot ((m_2 (x) - m_2 (y)) \times m_2 (y)) \de x \hspace{0.17em} \de
    y. \nonumber
  \end{align}
  Thus, adding and subtracting the quantity $(m_1 (x) - m_1 (y)) \times m_2
  (y)$, we infer 
  \begin{align}
    \left| \mathcal{H}_{\ep} (m_1) -\mathcal{H}_{\ep} (m_2) \right| &
    \leqslant  \iint_{\Omega \times \Omega} \frac{| \nu_{\ep} (x - y) |}{|x -
    y|} | m_1 (x) - m_1 (y) |  | m_1 (y) - m_2 (y)  | \de x \hspace{0.17em}
    \de y \nonumber\\
    &   + \iint_{\Omega \times \Omega} \frac{| \nu_{\ep} (x - y) |}{|x - y|}
    | (m_1 (x) - m_1 (y)) - (m_2 (x) - m_2 (y)) | | m_2 (y) | \de x
    \hspace{0.17em} \de y. \nonumber\\
    & \backassign  I_1 + I_2 .  \label{eq:estfordensity}
  \end{align}
  Let now $\tilde{m}_1, \tilde{m}_2 \in H^1 (\R^3 ; \R^3)$ be extensions of
  $m_1, m_2$ from $\Omega$ to $\R^3$. For the first integral $I_1$ in
  \eqref{eq:estfordensity} we have
  \begin{align}
    I_1 & \leqslant  \int_{\R^3} \frac{| \nu_{\ep} (h) |}{|h|} \left(
    \int_{\R^3} | \tilde{m}_1 (x + h) - \tilde{m}_1 (x) |  | \tilde{m}_1 (x +
    h) - \tilde{m}_2 (x + h) | \hspace{0.17em} \de x \right) \de h \nonumber\\
    & \leqslant  \int_{\R^3} | \nu_{\ep} (h) | \hspace{0.17em} \frac{\lVert
    \tau_h  \tilde{m}_1 - \tilde{m}_1 \rVert_{L^2 (\R^3)}}{|h|} 
    \hspace{0.17em} \| \tilde{m}_1 - \tilde{m}_2 \|_{L^2 (\R^3)}
    \hspace{0.17em} \de h \nonumber\\
    & \leqslant  \| \nabla \tilde{m}_1 \|_{L^2 (\R^3)}  \lVert \tilde{m}_1 -
    \tilde{m}_2 \rVert_{L^2 (\R^3)} \nonumber\\
    & \leqslant  C \hspace{0.17em} \|m_1 \|_{H^1 (\Omega)}  \|m_1 - m_2
    \|_{L^2 (\Omega)},  \label{eq:estfordensity01}
  \end{align}
  for some constant $C = C (\Omega) > 0$, where we have applied H{\"o}lder
  inequality and classical properties of Sobolev spaces,
  cf.~{\cite[Prop.~9.3]{brezis2011functional}}. Analogously, for the second
  integral $I_2$ in \eqref{eq:estfordensity} we find
  \begin{align}
    I_2 & \leqslant  \int_{\R^3} | \nu_{\ep} (h) |  \int_{\R^3} | \tilde{m}_2
    (x + h) | \frac{| (\tilde{m}_1 - \tilde{m}_2) (x) - (\tilde{m}_1 -
    \tilde{m}_2) (x + h) |}{|h|} \hspace{0.17em} \de x \hspace{0.17em} \de h
    \nonumber\\
    & \leqslant  \int_{\R^3} \left| \nu_{\ep} (h) \right| \lVert \tilde{m}_2
    \rVert_{L^2 (\R^3)} \frac{\lVert \tau_h (\tilde{m}_1 - \tilde{m}_2) -
    (\tilde{m}_1 - \tilde{m}_2) \rVert_{L^2 (\R^3)}}{|h|} \hspace{0.17em} \de
    h \nonumber\\
    & \leqslant  \lVert \nabla (\tilde{m}_1 - \tilde{m}_2) \rVert_{L^2
    (\R^3)} \lVert \tilde{m}_2 \rVert_{L^2 (\R^3)} \nonumber\\
    & \leqslant  C \hspace{0.17em} \lVert m_1 - m_2 \rVert_{H^1 (\Omega)}
    \lVert m_2 \rVert_{L^2 (\Omega)},  \label{eq:estfordensity02}
  \end{align}
  for some constant $C = C (\Omega) > 0$.
  
  To sum up, estimating \eqref{eq:estfordensity} through
  \eqref{eq:estfordensity01} and \eqref{eq:estfordensity02}, we obtain that,
  for every $m_1, m_2 \in H^1 (\Omega ; \R^3)$, there holds
  \begin{equation}
    \left| \mathcal{H}_{\ep} (m_1) -\mathcal{H}_{\ep} (m_2) \right| \leqslant
    C \hspace{0.17em} (\|m_1 \|_{H^1 (\Omega)} +\|m_2 \|_{H^1 (\Omega)}) 
    \|m_1 - m_2 \|_{H^1 (\Omega)} . \label{eq:equibound}
  \end{equation}
  Now, by density, for every $m \in H^1 (\Omega ; \R^3)$ there exists a
  sequence $(m_k)_{k \in \mathbb{N}}$ in $C^2_c  (\R^3 ; \R^3)$ such that
  \begin{equation}
    \lim_{k \to \infty} \|m - m_k \|_{H^1 (\Omega)} = 0. \label{eq:density}
  \end{equation}
  Also, by {\tmname{Step~1}}, we know that for every $m_k \in C^2_c  (\R^3 ;
  \R^3)$ one has
  \begin{equation}
    \lim_{\ep \to 0} \mathcal{H}_{\ep} (m_k) = \sum_{i = 1}^3 \int_{\Omega}
    m_k (x) \cdot (d_i \times \partial_i m_k (x)) \hspace{0.17em} \de x
    =\mathcal{H} (m_k) \label{eq:smoothcase}
  \end{equation}
  Moreover, in view of \eqref{eq:equibound}, for every $k \in \mathbb{N}$,
  there holds
  \begin{equation}
    \limsup_{\ep \to 0} \hspace{0.17em} \left| \mathcal{H}_{\ep} (m)
    -\mathcal{H}_{\ep} (m_k) \right| \leqslant C (\|m\|_{H^1 (\Omega)} +\|m_k
    \|_{H^1 (\Omega)})  \|m - m_k \|_{H^1 (\Omega)} . \label{eq:bound}
  \end{equation}
  Letting now $k \to + \infty$, \UUU by \eqref{eq:density} \EEE
  and by \eqref{eq:bound},
  \begin{equation}
    \lim_{k \to + \infty} \limsup_{\ep \to 0} \hspace{0.17em} \left|
    \mathcal{H}_{\ep} (m) -\mathcal{H}_{\ep} (m_k) \right| = 0,
    \label{eq:zero}
  \end{equation}
  thanks to the density \eqref{eq:density} and the fact that $\|m_k \|_{H^1
  (\Omega)}$ is bounded in $H^1$.
  
  \UUU In view of \eqref{eq:density}, \EEE we also have
  that $d_i \times \partial_i m_k \to d_i \times \partial_i m$ in $L^2 (\Omega
  ; \R^3)$ for $i = 1, 2, 3$. As a consequence, we have that
  \begin{equation}
    \lim_{k \to \infty} \mathcal{H} (m_k) =\mathcal{H} (m), \label{eq:cont}
  \end{equation}
  that is, the limit functional $\mathcal{H}$ is continuous with respect to
  the $H^1$-convergence.
  
  To sum up, by \eqref{eq:smoothcase}, \eqref{eq:zero} and \eqref{eq:cont},
  we obtain that
  \begin{align}
    \limsup_{\ep \to 0} \hspace{0.17em} \left| \mathcal{H}_{\ep} (m)
    -\mathcal{H}(m) \right| & \leqslant  \limsup_{\ep \to 0} \hspace{0.17em}
    \left( \left| \mathcal{H}_{\ep} (m) -\mathcal{H}_{\ep} (m_k) \right| +
    \hspace{0.17em} \left| \mathcal{H}_{\ep} (m_k) -\mathcal{H}(m_k) \right|
    \right) + \UUU | \mathcal{H}(m_k)
    -\mathcal{H}(m) |\EEE\nonumber\\
    & =  \left( \limsup_{\ep \to 0} \hspace{0.17em} \left| \mathcal{H}_{\ep}
    (m) -\mathcal{H}_{\ep} (m_k) \right| \right) + | \mathcal{H}(m_k)
    -\mathcal{H}(m) | . 
  \end{align}
  Passing to the limit for $k \rightarrow \infty$, we conclude that for
  every $m \in H^1 (\Omega ; \R^3)$, \UUU there holds \EEE $\mathcal{H}_{\ep} (m) \rightarrow
  \mathcal{H} (m)$ \UUU as \EEE $\varepsilon \rightarrow 0$, i.e., 
  \begin{equation*}
    \lim_{\ep \to 0} \iint_{\Omega \times \Omega} \frac{\nu_{\ep}  (x - y)}{|x
    - y|} \cdot (m (x) \times m (y)) \de x \hspace{0.17em} \de y = \sum_{i =
    1}^3 \int_{\Omega} m (x) \cdot (d_i \times \partial_i m (x)) \de x.
  \end{equation*}
  This completes the proof.
\end{proof}

\section{Compactness and \texorpdfstring{$\Gamma$}{Gamma}-convergence (proof of
Theorem~\ref{thm:gamma-conv})}\label{sec:Gammaconv}

\begin{proof}[\tmname{\bf Proof of Theorem~\ref{thm:gamma-conv}.{\tmem{i}},
{\tmname{\opt Compactness\cpt}}}]
  By the proof of Theorem \ref{thm:charact}
  ({\tmabbr{cf.~\eqref{eq:funestfepshepsineeps})}} we know that for every
  $\varepsilon > 0$ there holds
  \begin{equation}
    \frac{1}{2} \mathcal{F}_{\ep} (m_{\varepsilon}) \leqslant C \left(
    \mathcal{E}_{\ep} (m_{\varepsilon}) +\|m_{\varepsilon} \|_{L^2 (\Omega)}^2
    \right), \label{ineq:coerc}
  \end{equation}
  with $C > 0$ depending only on the homonymous constant in
  \ref{A1}. Since $(m_{\ep})_{\ep}$ consists of
  $\Stwo^2$-valued vector fields, i.e., $|m_{\ep} | = 1$ a.e. in $\Omega$ for
  every $\ep > 0$, from the assumption \eqref{eq:assumptionsupEeps} and
  \eqref{ineq:coerc} we infer that
  \begin{equation}
    \liminf_{\ep \rightarrow 0} \mathcal{F}_{\ep} (m_{\varepsilon}) \leqslant
    C, \label{eq:liminfestimatetemp}
  \end{equation}
  for some $C > 0$ depending only on $\liminf_{\ep \rightarrow 0} 
  \hspace{0.17em} \mathcal{E}_{\ep} (m_{\ep})$, $\Omega$, and on the
  homonymous constant in \ref{A1}. 
  \UUU In particular, in view of assumption \ref{G4}, recalling the definition of the kernels $(\rho_\ep^{\rm rad})_\ep$ in \eqref{eq:radial-k} and setting $$\tilde{\rho}_\ep(x):=\frac{\rho_\ep^{\rm rad}(|x|)}{\|\rho_\ep^{\rm rad}(|\cdot|)\|_{L^1(\R^3)}},$$ from estimate
  \eqref{eq:liminfestimatetemp} we infer
  \begin{equation}
      \label{eq:liminfestimatetemp2}
      \liminf_{\ep\to 0}\iint_{\Omega \times \Omega} \tilde{\rho}_{\ep}  (x - y)  \frac{|m (x) - m (y)
    |^2}{|x - y|^2} \hspace{0.17em} \de x \hspace{0.17em} \de y\leqslant \frac{C}{\kappa}.
  \end{equation}
  \EEE
  
  Estimate
  \eqref{eq:liminfestimatetemp2} allows us to invoke the compactness result for
  the symmetric exchange energy  established
  in~{\cite[Theorem 1.2]{ponce2004estimate}} for the case of radial kernels, which implies that if
  $(m_{\varepsilon})_{\varepsilon}$ is a bounded family in $L^2 (\Omega ;
  \R^3)$ such that $\mathcal{F}_{\ep} (m_{\varepsilon}) \leqslant C$ for some
  positive constant $C$, then there exists $m \in H^1 (\Omega ; \R^3)$ and a
  subsequence $\ep_j \to 0$ such that $m_{\ep_j} \to m$ strongly in $L^2
  (\Omega ; \R^3)$. In particular, we obtain that (up to a further
  subsequence) $m_{\ep_j} \to m$ a.e. as $j \to \infty$ and, therefore, also
  the limiting function satisfies the constraint, that is, $m \in L^2 (\Omega ;
  \Stwo^2)$.
  \end{proof}

\begin{proof}[\tmname{\bf Proof of Theorem~\textup{\ref{thm:gamma-conv}}.{\tmem{ii}}, ($\Gamma$-{\tmname{convergence}})}]
  We split the proof into two classical steps. First, we derive the
  $\Gamma$-$\mathrm{limsup}$ inequality, then we focus on the
  $\Gamma$-$\mathrm{liminf}$ inequality.{\medskip}
  
  {\noindent}{\tmname{Limsup inequality.}} It is sufficient to note that the
  constant family $ m_{\ep}:=m$ for every $\ep$, is a recovery sequence. Indeed, on the one hand, if $m
  \in L^2 (\Omega ; \Stwo^2) \setminus H^1 (\Omega ; \Stwo^2)$ then
  Theorem~\ref{thm:charact} gives that $\mathcal{E}_{\ep} (m) \rightarrow +
  \infty$; on the other hand, if $m \in H^1 (\Omega ; \Stwo^2)$ then by
  Corollary~\ref{cor:point-conv}, we conclude that
  \[ \lim_{\varepsilon \to 0} \mathcal{E}_{\varepsilon} (m) =\mathcal{E} (m) .
  \]
  {\medskip}{\noindent}{\tmname{Liminf inequality.}} Let $(m_{\ep})_{\ep}
  \subset L^2 (\Omega ; \Stwo^2)$ be such that $m_{\ep} \to m$ strongly in
  $L^2 (\Omega ; \Stwo^2)$ as $\ep \to 0$. If $\liminf_{\varepsilon
  \rightarrow 0} \mathcal{E}_{\varepsilon} (m_{\varepsilon}) = + \infty$,
  there is nothing to prove. If $\liminf_{\varepsilon \rightarrow 0}
  \mathcal{E}_{\varepsilon} (m_{\varepsilon}) < + \infty$ then, by
  \eqref{ineq:coerc}, we know that
  \begin{equation}
    \liminf_{\ep \to 0} \mathcal{F}_{\ep} (m_{\varepsilon}) \leqslant C \left(
    1 + \liminf_{\ep \to 0} \mathcal{E}_{\ep} (m_{\varepsilon}) \right) < +
    \infty, \label{ineq:coerc2}
  \end{equation}
  for some constant $C$ depending only on $\Omega$ and the homonymous constant
  in \ref{A1}. Also, by the compactness result
  (Theorem~\ref{thm:gamma-conv}.{\tmem{i}}), $m \in H^1 \left( \Omega ;
  \Stwo^2 \right)$, and from {\cite[Lemma 8]{ponce2004new}} used in the case
  $\omega (t) = t^2$, we get the existence of a finite Radon measure $\mu \in
  \mathcal{M} (\Stwo^2)$, with $\mu \left( \Stwo^2 \right) = 1$, such that,
  \UUU possibly \EEE up to \UUU the extraction of a non-relabelled \EEE  subsequence, there holds
  \begin{equation}
    \label{Sym:liminf} \liminf_{\ep \to 0} \mathcal{F}_{\ep} (m_{\varepsilon})
    \geqslant \int_{\Omega} \left( \int_{\Stwo^2} | \partial_{\sigma} m (x)
    |^2 \de \mu (\sigma) \right) \de x.
  \end{equation}
  Therefore, it only remains to show the analogous $\mathrm{liminf}$
  inequality for the nonlocal antisymmetric exchange term, i.e., that
  \begin{equation}
    \liminf_{\ep \to 0} \mathcal{H}_{\varepsilon} (m_{\varepsilon}) \geqslant
    \mathcal{H} (m) . \label{eq:ourgoalliminf}
  \end{equation}
  We show something more, as stated in the next \UUU claim. 
  
  \noindent\emph{Claim}:
    \begin{equation}
      \lim_{\varepsilon \rightarrow 0} \mathcal{H}_{\varepsilon}
      (m_{\varepsilon}) =\mathcal{H} (m) . \label{eq:ourgoalliminbetter}
    \end{equation}\EEE

  \UUU In order to prove our claim, \EEE we introduce a family of mollifiers $(\eta_k)_{k \in \mathbb{N}}$ in
  $C^{\infty}_c (\R^3)$ build in the usual way: we consider an element $\eta
  \in C^{\infty}_c (\R^3)$, $\eta \geqslant 0$, $\| \eta \|_{L^1 (\R^3)} = 1$,
  $\mathrm{supp} \eta \subset B_1 (0)$, and we set $\eta_k (x) \assign k^3
  \eta (k x)$. Also, for every sufficiently small $\delta > 0$ we set
  $\Omega_{\delta} \assign \{x \in \Omega \of \mathrm{dist} (x, \partial
  \Omega) > \delta\}$. For every $k > 1 / \delta$ the regularization of $m_{\ep}$ given by $m_{\ep, k}
  \assign m_{\ep} \ast \eta_k$ is well-defined on
  $\Omega_{\delta}$. By classical properties of mollifiers, we know
  that for a given $\ep > 0$ there holds
  \begin{align}
    m_{\ep, k} \assign m_{\ep} \ast \eta_k & \xrightarrow{k \rightarrow +
    \infty}  m_{\ep} \quad \text{strongly in } L^2 (\Omega_{\delta} ; \R^3), 
    \label{moll:prop}\\
    m_k \assign m \ast \eta_k & \xrightarrow{k \rightarrow + \infty} 
    m_{\phantom{\varepsilon}} \quad \text{strongly in } H^1 (\Omega_{\delta} ;
    \R^3) .  \label{moll:prop2}
  \end{align}
  In writing the previous relations, we agreed (with a safe abuse of notation)
  that the $\Stwo^2$-valued vector fields $m$ and $m_{\varepsilon}$ are
  extended by zero outside of $\Omega$.
  
  In what follows, it will be convenient to denote by
  $\mathcal{H}_{\varepsilon, \delta}$ and $\mathcal{H}_{\delta}$ the analog of
  the functionals $\mathcal{H}_{\varepsilon}$ and $\mathcal{H}$ when defined
  on space of functions defined on the compact set $\bar{\Omega}_{\delta}
  \subset \Omega$ instead of $\Omega$. Thus, for example,
  \begin{equation*}
    \mathcal{H}_{\ep, \delta} (m_{\ep, k}) = \iint_{\UUU\bar{\Omega}_{\delta}^2}
    \frac{\nu_{\ep}  (x - y)}{|x - y|} m_{\ep, k} (x) \times m_{\ep, k} (y)
    \de x \de y.
  \end{equation*}
  For the nonlocal antisymmetric energy $\mathcal{H}_{\varepsilon, \delta}$ we
  then observe that the uniform convergence result stated in
  Theorem~\ref{thm:point-conv} assures that
  \begin{equation}
    \lim_{\ep \to 0} \mathcal{H}_{\ep, \delta} (m_{\ep, k})
    =\mathcal{H}_{\delta} (m_k) = \sum_{i = 1}^3 \int_{\Omega_{\delta}} m_k
    (x) \cdot (d_i \times \partial_i m_k (x)) \de x. \label{liminf:regular}
  \end{equation}
  Indeed, since $m_{\ep} \to m$ in $L^2 (\Omega ; \Stwo^2)$ as $\ep \to 0$, we
  have that for every $k \in \mathbb{N}$ there holds $m_{\ep, k} = m_{\ep}
  \ast \eta_k \xrightarrow{\varepsilon \rightarrow 0} m_k = m \ast \eta_k$ in
  $C^2 (\bar{\Omega}_{\delta} ; \R^3)$, and therefore, we are in the
  hypotheses of \eqref{eq:Hepspointwiseconvunif} in
  Theorem~\ref{thm:point-conv} applied to $\bar{\Omega}_{\delta}$.
  
  Note that, by \eqref{moll:prop2} we have
  $\mathcal{H} (m) = \lim_{\delta \rightarrow 0} (\lim_{k \rightarrow +
  \infty} \mathcal{H}_{\delta} (m_k))$ and, therefore, our goal
  \eqref{eq:ourgoalliminbetter} boils down to showing the equality
  \begin{equation}
    \lim_{\varepsilon \rightarrow 0} \mathcal{H}_{\varepsilon}
    (m_{\varepsilon}) = \lim_{\delta \rightarrow 0} \lim_{k \rightarrow +
    \infty} \lim_{\ep \to 0} \mathcal{H}_{\ep, \delta} (m_{\ep, k}),
    \label{eq:ourgoalliminfnew}
  \end{equation}
  because, by \eqref{liminf:regular}, the RHS of \eqref{eq:ourgoalliminfnew}
  is nothing but $\mathcal{H} (m)$. To show equality
  \eqref{eq:ourgoalliminfnew} we split $\mathcal{H}_{\varepsilon} \left(
  m_{\ep} \right)$ under the form
  \begin{equation}
    \mathcal{H}_{\varepsilon} (m_{\varepsilon}) =\mathcal{H}_{\varepsilon,
    \delta} \left( m_{\ep, k} \right) + (\mathcal{H}_{\varepsilon}
    (m_{\varepsilon}) -\mathcal{H}_{\delta, \varepsilon} (m_{\varepsilon})) +
    \left( \mathcal{H}_{\delta, \varepsilon} (m_{\varepsilon})
    -\mathcal{H}_{\varepsilon, \delta} \left( m_{\ep, k} \right) \right)
    \label{eq:newgoal}
  \end{equation}
  and then show that
  \begin{align}
    \lim_{\delta \rightarrow 0} \lim_{k \rightarrow + \infty}
    \limsup_{\varepsilon \rightarrow 0} | \mathcal{H}_{\varepsilon}
    (m_{\varepsilon}) -\mathcal{H}_{\delta, \varepsilon} (m_{\varepsilon}) | &
    =  0,  \label{eqs:usefullimitrels1}\\
    \lim_{\delta \rightarrow 0} \lim_{k \rightarrow + \infty}
    \limsup_{\varepsilon \rightarrow 0} \left| \mathcal{H}_{\delta,
    \varepsilon} (m_{\varepsilon}) -\mathcal{H}_{\varepsilon, \delta} \left(
    m_{\ep, k} \right) \right| & =  0.  \label{eqs:usefullimitrels2}
  \end{align}
  Indeed, as soon as we achieve
  \eqref{eqs:usefullimitrels1}--\eqref{eqs:usefullimitrels2}, passing to the
  limit $\lim_{\delta \rightarrow 0} \lim_{k \rightarrow + \infty}
  \limsup_{\varepsilon \rightarrow 0}$ in \eqref{eq:newgoal} we
  conclude.{\smallskip}
  
  {\smallskip}{\noindent}{\tmname{Step 1. Proof of
  \eqref{eqs:usefullimitrels1}}}. \ Using the fact that the kernel
  $\nu_{\varepsilon}$ is odd, we \UUU obtain \EEE that
  \begin{align}
    \mathcal{H}_{\varepsilon} (m_{\varepsilon}) -\mathcal{H}_{\delta,
    \varepsilon} (m_{\varepsilon}) & =   \iint_{\Omega^2\setminus \Omega_{\delta}^2}  \frac{\nu_{\ep}  (x - y)}{|x - y|} \cdot
    \left( m_{\ep} (x) \times m_{\ep} (y) \right)  \de x \de y \nonumber
  \end{align}
  and, therefore,
  \begin{align}
   | \mathcal{H}_{\varepsilon} (m_{\varepsilon}) -\mathcal{H}_{\delta,
     \varepsilon} (m_{\varepsilon}) | & \leqslant 3 \left( C \iint_{\R \times
     (\Omega \backslash \Omega_{\delta})} \left| \nu_{\ep} (x - y) \right| 
     \de x \de y \right)^{1 / 2} \cdot \left( \mathcal{F}_{\ep}
     (m_{\varepsilon}) \right)^{1 / 2} \notag \\
     &= 3 C^{1 / 2} \left( \mathcal{F}_{\ep}
     (m_{\varepsilon}) \right)^{1 / 2} | \Omega \backslash \Omega_{\delta}
     |^{1 / 2} \notag.
     \end{align}
  Passing to the limit in the previous inequality, we \UUU infer \EEE
  \eqref{eqs:usefullimitrels1}.{\smallskip}
  
  {\noindent}{\tmname{Step 2. Proof of \eqref{eqs:usefullimitrels2}}}. We \UUU observe \EEE
  that
  \begin{align}
    \mathcal{H}_{\delta, \varepsilon} (m_{\varepsilon})
    & -\mathcal{H}_{\varepsilon, \delta} \left( m_{\ep, k} \right)\notag \\
     & = 
    \iint_{\Omega_{\delta}^2} \frac{\nu_{\ep}  (x - y)}{|x - y|} \cdot \left(
    m_{\ep} (x) \times m_{\ep} (y) - m_{\ep, k} (x) \times m_{\ep, k} (y)
    \right)  \de x \de y \nonumber\\
    & =  \iint_{\R^3 \times \R^3} \left( \iint_{\Omega_{\delta}^2}
    \frac{\nu_{\ep} (x - y)}{|x - y|} \cdot \left( m_{\ep} (x) \times m_{\ep}
    (y) \right)  \de x \de y \right) \eta_k (z) \eta_k (w) \de z
    \hspace{0.17em} \de w \nonumber\\
    &   - \iint_{\R^3 \times \R^3} \left( \iint_{(- z + \Omega_{\delta})
    \times (- w + \Omega_{\delta})} \frac{\nu_{\ep} (x - y)}{|x - y|} \cdot
    \left( m_{\ep} (x) \times m_{\ep} (y) \right)  \de x \de y \right) \eta_k
    (z) \eta_k (w) \de z \hspace{0.17em} \de w \nonumber\\
    & =  \iint_{\R^3 \times \R^3} \left[ \left( \iint_{\Omega^2}
    \frac{\nu_{\ep} (x - y)}{|x - y|} \cdot \left( m_{\ep} (x) \times m_{\ep}
    (y) \right) \chi_{(z, w)} (x, y)  \de x \de y \right) \right] \eta_k (z)
    \eta_k (w) \de z \hspace{0.17em} \de w,  \label{eq:step2toestimate}
  \end{align}
  where we set $\chi_{(z, w)} (x, y) \assign 1_{\Omega_{\delta}} (x)
  1_{\Omega_{\delta}} (y) - 1_{(- z + \Omega_{\delta})} (x) 1_{(- w +
  \Omega_{\delta})} (y)$. Before continuing estimating
  \eqref{eq:step2toestimate}, we need the following inequality, whose proof is
  postponed to {\tmname{Step~3}} below:
  \begin{align}
    \iint_{\R^3 \times \R^3} 
    & \left( \iint_{\Omega^2}
    \left| \nu_{\ep} (x - y) \right| | \chi_{(z, w)} (x, y) |  \de x \de y
    \right) \eta_k (z) \eta_k (w) \de z \de w \notag\\
    & \qquad\qquad\qquad\qquad \leqslant  2 \int_{\R^3}
    \int_{\Omega} | 1_{\Omega_{\delta}} (x) - 1_{(- z + \Omega_{\delta})} (x)
    |  \de x \eta_k (z) \de z.  \label{eq:tobeusedsoon}
  \end{align}
  Having \eqref{eq:tobeusedsoon} at our disposal, \UUU by combining H{\"o}lder inequality, \ref{A1}, as well as
  \eqref{eq:step2toestimate}, we deduce \EEE
\begin{align}
       | \mathcal{H}_{\delta, \varepsilon} (m_{\varepsilon})
       & -\mathcal{H}_{\varepsilon, \delta} ( m_{\ep, k}) |\notag \\
      \leqslant &  \left( 2\mathcal{F}_{\ep} (m_{\varepsilon}) \right)^{1 / 2}
       \iint_{\left( \R^3 \right)^2} \left( \iint_{\Omega^2} \frac{\left|
       \nu_{\ep} (x - y) \right|^2}{\rho_{\varepsilon} (x - y)} | \chi_{(z,
       w)} (x, y) |  \de x \de y \right)^{1 / 2} \eta_k (z) \eta_k (w) \de z
       \hspace{0.17em} \de w \nonumber\\
       \leqslant & C^{1 / 2} \left( 2\mathcal{F}_{\ep} (m_{\varepsilon})
       \right)^{1 / 2} \iint_{\left( \R^3 \right)^2} \left( \iint_{\Omega^2}
       \left| \nu_{\ep} (x - y) \right| | \chi_{(z, w)} (x, y) |  \de x \de y
       \right)^{1 / 2} \eta_k (z) \eta_k (w) \de z \de w \nonumber\\
       \leqslant & C^{1 / 2} \left( 2\mathcal{F}_{\ep} (m_{\varepsilon})
       \right)^{1 / 2} \left( \iint_{\left( \R^3 \right)^2} \left(
       \iint_{\Omega^2} \left| \nu_{\ep} (x - y) \right| | \chi_{(z, w)} (x,
       y) |  \de x \de y \right) \eta_k (z) \eta_k (w) \de z \hspace{0.17em}
       \de w \right)^{1 / 2} \nonumber\\
       \leqslant & 2 C^{1 / 2} \left( \mathcal{F}_{\ep} (m_{\varepsilon})
       \right)^{1 / 2} \left( \int_{\R^3} \int_{\Omega} | 1_{\Omega_{\delta}}
       (x) - 1_{(- z + \Omega_{\delta})} (x) | \eta_k (z) \de x  \de z
       \right)^{1 / 2} \nonumber\\
       = & 2 C^{1 / 2} \left( \mathcal{F}_{\ep} (m_{\varepsilon}) \right)^{1 /
       2} \left( \int_{B_1} \| 1_{\Omega_{\delta}} (\cdot) -
       1_{\Omega_{\delta}} (\cdot - z / k) \|_{L^1 (\Omega)}  \de z \right)^{1
       / 2} \xrightarrow{k \rightarrow \infty} 0,  \label{eq:finalresult}
\end{align}
  the convergence of \eqref{eq:finalresult} to $0$ being a consequence of
  \UUU Lebesgue's Dominated Convergence Theorem, as well as of the continuity of the $L^2$-norm with respect to translations.{\medskip}
  
  {\noindent}{\tmname{Step~ 3}}. {\tmname{Proof of \eqref{eq:tobeusedsoon}}}.
  \UUU We first \EEE rewrite the function $\chi_{(z, w)} (x, y)$ \UUU as \EEE
  \[ \chi_{(z, w)} (x, y) = [1_{\Omega_{\delta}} (x) - 1_{(- z +
     \Omega_{\delta})} (x)] 1_{\Omega_{\delta}} (y) + 1_{(- z +
     \Omega_{\delta})} (x) [1_{\Omega_{\delta}} (y) - 1_{(- w +
     \Omega_{\delta})} (y)], \]
  from which we \UUU infer \EEE the following estimate:
  \[ | \chi_{(z, w)} (x, y) | \leqslant | 1_{\Omega_{\delta}} (x) - 1_{(- z +
     \Omega_{\delta})} (x) | + | 1_{\Omega_{\delta}} (y) - 1_{(- w +
     \Omega_{\delta})} (y) | . \]
  \UUU In particular, \EEE
  \begin{align}
    \iint_{\Omega^2} \left| \nu_{\ep} (x - y) \right| | \chi_{(z, w)} (x, y) |
    \de x \de y \leqslant & \iint_{\Omega \times \R^3} \left| \nu_{\ep} (x
    - y) \right| | 1_{\Omega_{\delta}} (x) - 1_{(- z + \Omega_{\delta})} (x) |
    \de x \de y \nonumber\\
    &  \qquad + \iint_{\R^3 \times \Omega} \left| \nu_{\ep} (x - y) \right|
    | 1_{\Omega_{\delta}} (y) - 1_{(- w + \Omega_{\delta})} (y) | \de x \de y
    \nonumber\\
    & =  \int_{\Omega} | 1_{\Omega_{\delta}} (x) - 1_{(- z +
    \Omega_{\delta})} (x) | + | 1_{\Omega_{\delta}} (x) - 1_{(- w +
    \Omega_{\delta})} (x) |  \de x. \nonumber
  \end{align}
  The previous inequality, once integrated on $\R^3 \times \R^3$ with respect
  to the measure $\eta_k (z) \eta_k (w) \de z \hspace{0.17em} \de w$ \UUU yields \EEE
  \eqref{eq:tobeusedsoon}.
\end{proof}

\section*{Acknowledgements}
\noindent{\sc G.DiF.} and {\sc R.G.} acknowledge support from the Austrian Science Fund (FWF) through the project {\emph{Analysis and Modeling of Magnetic
    Skyrmions}} (grant 10.55776/P34609). The research of {\sc E.D.} has been supported by the Austrian Science Fund (FWF) through grants 10.55776/F65, 10.55776/V662,
10.55776/P35359 and 10.55776/Y1292, as well as from BMBWF through the OeAD/WTZ project CZ 09/2023. 
    {\sc G.DiF.}
thanks TU Wien and MedUni Wien for their support and hospitality.

\bibliographystyle{abbrv}
\bibliography{DDFG24}
\end{document}